# Two-time, response-excitation moment equations for a cubic half-oscillator under Gaussian and cubic-Gaussian colored excitation.
# Part 1: The monostable case

by
G. A. Athanassoulis([1]),    I. C. Tsantili    and    Z. G. Kapelonis
*mathan@central.ntua.gr, ivits@central.ntua.gr, zkapel@central.ntua.gr*
*National Technical University of Athens, School of Naval Architecture and Marine Engineering,*
*Section of Ship and Marine Hydrodynamics, Heroon Polytechniou 9, 157 73, Athens, Greece*

**Abstract**

In this paper a new method is presented for the formulation and solution of two-time, response-excitation moment equations for a nonlinear dynamical system excited by colored, Gaussian or non-Gaussian processes. Starting from equations for the two-time moments (e.g. for $C_{xy}(t,s)$, $C_{xx}(t,s)$), the method uses an exact time-closure condition, in addition to a Gaussian moment closure, in order to obtain a closed, non-local in time (causal) subsystem for the one-time ($t=s$) moments. After solving this causal system, the two-time moments can be calculated for all $(t,s)$ pairs as well. The present method differs essentially from the classical Itô/FPK approach since it does not involve any specific assumptions regarding the correlation structure of the excitation. In the case where the input random process can be obtained as the solution of an Itô equation (as, e.g., happens with an Ornstein-Uhlenbeck process), the proposed non-local system is localized, leading to moment equations identical with the usual ones. The closed, non-local in time, moment system is numerically solved by means of an appropriate, two-scale, iterative scheme, and numerical results are presented for two families of colored stochastic excitations. The results are confirmed by means of extensive Monte Carlo simulations. It is found that both the correlation time and the details of the shape of the input random function affect appreciable the response covariance. In the present paper we focus on a monostable cubic half-oscillator, excited by a smoothly-correlated, linear-plus-cubic-Gaussian (non-Gaussian) random input. The bistable case, as well as more general nonlinear systems can be treated by the same method, provided that a more elaborate moment closure will be used.

---

**Highlights:**
- We introduce two-time, response-excitation moment equations for nonlinear random equations
- The excitation process is smoothly correlated (colored noise), with arbitrary correlation structure
- The moment equations form a non-local in time, causal system
- The moment system is solved numerically by a two-scale, iterative approach
- Numerical results are confirmed by Monte Carlo simulations
- It is found that both correlation time and the details of the shape of the input random function affect appreciable the response covariance

# Two-time, response-excitation moment equations for a cubic half-oscillator under Gaussian and cubic-Gaussian colored excitation.
# Part 1: The monostable case

## Table of contents

# 1. Introduction

In many engineering problems the input data are not (or cannot be) accurately known, introducing a degree of uncertainty in the analysis of the examined system. In such cases probabilistic modeling is usually invoked in order to quantify the uncertainty. Some or all of the data elements are modeled as random variables or random functions, in accordance with their nature, and the response is probabilistically predicted as well. Input (data) random functions may appear as external or parametric excitations in the dynamical equations governing the evolution of the examined system.

The probabilistic characterization of data random functions is an important and difficult procedure, which may take years of intensive research in order to reach a satisfactory level. Sea-wave loads, wind loads, and seismic ground motion are three well-known examples. The complete probabilistic characterization of a random function encodes a great amount of information. It requires the knowledge of the hierarchy of the probability distributions of all orders or, equivalently, the knowledge of the characteristic functional [1]. Usually only partial knowledge of a data random function, say $y(t;\theta)$, is available, either in the form of some moments, e.g.

$$m_y(t) = \mathrm{E}^\theta[y(t;\theta)], \quad C_{yy}(t,s) = \mathrm{E}^\theta[(y(t;\theta) - m_y(t)) \cdot (y(s;\theta) - m_y(s))],$$

or in the form of some low-order probability density functions (pdfs), e.g. $f_{y(t)}(\beta)$, $f_{y(t)y(s)}(\beta_1, \beta_2)$. Here and subsequently $\theta$ is the stochastic argument and $\mathrm{E}^\theta[\bullet]$ is the mean value operator.

Assuming that, in a given problem, the probabilistic knowledge of the data has reached a widely acceptable status, the next step is to exploit it in order to quantify the uncertainty of the response. Various approaches have been developed for this purpose, usually characterized by means of the equations that are formulated and solved. The Fokker-Planck-Kolmogorov (FPK) equation approach (essentially equivalent to the Itô stochastic differential equations (SDEs) approach), the moment equations approach, the statistical linearization, and the Monte Carlo (MC) simulation are probably the favorite ones to the engineering community. Nowadays they are all well established and well described in various excellent books dealing with stochastic systems or random vibrations; see, e.g., [1], [2], [3], [4], [5], [6], [7], [8]. A restrictive feature of the Itô/FPK approach is that all input random functions are assumed to be Gaussian white noises. Generalizations permitting to consider other kind of noises (e.g., Poisson noise) have been developed as well; see, e.g., [1], [9], [10], [11]. Although this is an important extension, it does not go beyond delta-correlated excitations, which is too restrictive for many engineering applications. The range of applicability of the Itô/FPK approach has been further extended by means of the stochastic averaging approach ([12], [13], [14], [15], [16]) and the filtering approach (state-space augmentation, [17], [18], [19], [20]), which permit to consider various classes of smoothly-correlated (colored) input random functions. The latter method, subsequently called the Itô/filtering approach, is able to treat colored data processes that can be obtained as the output of an Itô SDE (the shaping filter).

A desirable feature of any method of uncertainty quantification is to be able to derive a partial (ideally, best possible) probabilistic characterization of the response, fully exploiting the available information concerning the data random functions. This is feasible for MC simulations for various classes of Gaussian and non-Gaussian random data functions (see, e.g., [3], [21], [22], [23], [24], [25]), but not always possible (or easy) for the other approaches, except for specific cases. One such case, concerning moments, refers to linear systems under random external excitation. In this case, the knowledge of some moments of the excitation results in the determination of the same moments for the response ([26], [27], [28], [29]). Another case, concerning moments and pdfs, refers to quite general dynamical systems, linear or non-linear, excited by delta-correlated processes, or processes that can be modeled as the output of some filter excited by other delta-correlated processes. Such problems can be solved by means of the Itô/FPK or Itô/filtering approach, or the various generalizations of them. Besides, in this case, it is possible to systematically formulate moment equations up to any desirable order. Of course, because of the non-linearity, the (truncated) moment system is not closed and thus some closure scheme should be invoked. The simplest one is the Gaussian closure, introduced by Goodman and Whittle in 50's and extensively used thenceforth in the study of random vibrations and more general stochastic dynamical systems [28], [7]. It has been found that it works well for monostable oscillators, while for bistable ones may lead to inadequate or erroneous results [30], [31]. Also, many types of non-Gaussian closures have been devised and used for treating moment equations coming from nonlinear stochastic systems under delta-correlated excitation. Among them we mention the cumulant-neglect closure ([32], [4], [33], [28], [7]), closures based on an a priori selection of a specific parametric model for the underlying pdf, including the quasi-moment neglect method [1], [34], [35], [36], [37], [38], the information closure [39], [40], the polynomial-Gaussian closure [41], [42], and the "power $p$ of order $n$" closure, first introduced in statistical mechanics and turbulence under the name "Kirkwood superposition approximation" (see [43] for a relevant review), and recently systematized by Hausken and Moxens [44], who also coined the name "power $p$ of order $n$" closure.

In this paper we focus on the formulation of and the solution to the two-time, response-excitation (RE) moment equations for a nonlinear dynamical system excited by colored, Gaussian or non-Gaussian processes. More precisely, we formulate, elaborate and solve moment equations corresponding to a first-order scalar equation with a cubic non-linearity, excited by a smoothly-correlated, linear-plus-cubic-Gaussian (non-Gaussian) random function. The derivation of the moment equations is straightforward, based on repeated use of the dynamical equation multiplied by the excitation or the response function, and being averaged (generalizing an approach usually applied to linear systems [29]). Our method of formulating moment equations differs in two respects from the classical methods using the Itô/FPK/filtering approaches. First, we start from the two-time moment equations, e.g. equations for $C_{xy}(t,s)$, $C_{xx}(t,s)$, and then we obtain a non-local in time subsystem for the one-time ($t=s$) moments using time closure (see below), while in the Itô/FPK/filtering approach the moment problem is formulated directly with respect to one-time moments, leading to a local in time, nonlinear system of ordinary differential equations (ODEs). Second, we do not need any specific assumptions regarding the correlation structure of the excitation. Nevertheless, if the excitation considered fits in the Itô/filtering scheme (as, e.g., in the case of an Ornstein-Uhlenbeck (OU) input process) our non-local equations are localized leading to moment equations identical with the usual ones. Accordingly, our method appears as a natural generalization of the Itô/filtering approach. In the present paper we focus on the monostable cubic half-oscillator. The bistable case, as well as more general nonlinear systems will be treated elsewhere [45].

The two-time moment system is not closed in two ways. On the one hand, it contains higher-order moments, as also the usual, one-time moment systems do. On the other hand, the involved moments depend on two time variables, while time differentiation applies only to one of them. Accordingly, two closure procedures are needed in order to obtain a closed, solvable system: a moment closure and a time closure. In this work the simplest moment closure, the Gaussian one, has been considered, the emphasis being put on the implementation of the time closure. We show that an exact time closure is possible without any additional assumptions. The obtained closed moment system is nonlinear and non-local in time, belonging to a family of functional differential equations introduced by Volterra and Tonelli in the beginning of the twentieth century. Nowadays such type of equations are studied under the name of causal differential equations/systems [46], [47], or Volterra functional differential equations [48], [49]. The closed moment system is numerically solved by means of an appropriate, two-scale, iterative scheme, and numerical results are presented for two families of colored random excitations. Numerical results indicate that the output two-time moments are signi-ficantly affected both by the correlation time and the details of the shape of the input random function, which in general cannot be taken into account in the Itô/FPK/filtering approaches. The obtained results compare satisfactorily with extensive Monte Carlo simulations results to the extent allowed by the Gaussian closure.

The present work, apart from its significance as a new method of formulating and solving two-time moment equations, constitutes a building block of a new general approach to the probabilistic study of dynamical systems under general, colored (smoothly-correlated), random excitation, aiming at the calculation of the joint RE pdf. This approach initiated by Athanassoulis and Sapsis [50], [51], [52], and has been elaborated further by Venturi, Karniadakis *et al* [53], [54], who coined the name response-excitation theory, adopted by the present authors as well. In the context of the general response-excitation (RE) theory, moment equations of the type studied herewith are consi-dered locally in the RE phase space, serving as closure conditions for the differential constraint satisfied by the joint RE pdf. Preliminary results in this direction have been published recently [55], [56] and will be fully elaborated and presented in forthcoming papers.

The following **abbreviations** will be frequently used in the subsequent sections

| | |
|---|---|
| FPK | Fokker-Planck-Kolmogorov (equation) |
| Gf | Gaussian filter (a Gaussian random process) |
| MC | Monte Carlo (simulation) |

ODE(s)  ordinary differential equation(s)
OU      Ornsten-Uhlenbeck (process)
pdf(s)  probability density function(s)
RDE(s)  random differential equation(s), i.e., ordinary differential equation(s) driven by smoothly-correlated random functions
RE      response-excitation (theory, moments, phase space)
SDE(s)  stochastic differential equation(s), i.e., ordinary differential equation(s) driven by delta-correlated stochastic functions

**Notation**

$$m_{x_0} = \mathrm{E}^{\theta}[x_0(\theta)], \quad m_x(t) = \mathrm{E}^{\theta}[x(t;\theta)], \quad m_y(t) = \mathrm{E}^{\theta}[y(t;\theta)],$$

$$C^{j_1 j_2}_{x\,y}(t,s) = \mathrm{E}^{\theta}\left[(x(t;\theta) - m_x(t))^{j_1} \cdot (y(s;\theta) - m_y(s))^{j_2}\right], \quad C^{11}_{xy}(t,s) = C_{xy}(t,s),$$

$$C^{1\,k}_{x_0\,x}(t_0,t) = C^{k}_{x_0\,x}(t) = \mathrm{E}^{\theta}\left[(x_0(\theta) - m_{x_0}) \cdot (x(t;\theta) - m_x(t))^{k}\right], \quad C^{1}_{x_0\,x}(t) = C_{x_0\,x}(t).$$

Similar notation is also used for the non-central second-order moments (cross-correlations), $R^{j_1 j_2}_{x\,y}(t,s)$, $R^{k}_{x_0\,x}(t)$, etc. The above notation is unambiguous except for the case $x = y$, $t = s$. In this case the moment function $C^{j_1 j_2}_{x\,x}(t,t)$ can be written by various, apparently different but fully equivalent, ways; all couples $(j_1, j_2)$ with $j_1 + j_2 = j$ represent the same $j$-th order central moment of $x(t;\theta)$. For example,

$$C^{21}_{xx}(t,t) = C^{12}_{xx}(t,t) = C^{30}_{xx}(t,t) = C^{03}_{xx}(t,t) = \mathrm{E}^{\theta}\left[(x(t;\theta) - m_x(t))^{3}\right].$$

## 2. The cubic half-oscillator and the random excitation

The cubic half-oscillator excited by delta-correlated processes has been studied by many authors. Hasofer and Grigoriu [30], [31] studied the case of Gaussian white-noise excitation, solving the corresponding moment problem and commenting on the properties of the moment closures; Wojtkiewicz, Grigoriu *et al* [9], [22] generalized various techniques developed for the case of Gaussian excitation to systems driven by Poisson or Gaussian plus Poisson white-noises, and studied the cubic half oscillator as an example. The same problem, under OU excitation, has been treated by Jung, Risken and others [57], [58], in the context of the filtering approach (augmented state space), with emphasis on the bistable case. Furthermore, the same problem has been also studied using approximate one-dimensional FPK equation in conjunction with the short relaxation time approximation [59], [60].

In the present paper we study a cubic half-oscillator described by the RDE

$$\dot{x}(t;\theta) = \mu_1 \cdot x(t;\theta) + \mu_3 \cdot x^3(t;\theta) + \kappa_1 \cdot y(t;\theta) + \kappa_3 \cdot y^3(t;\theta), \tag{1a}$$

$$x(t_0;\theta) = x_0(\theta), \tag{1b}$$

where $y(t;\theta)$ is a given, smoothly-correlated, Gaussian random function, and the initial condition $x_0(\theta)$ is a given random variable, independent from the process $y(t;\theta)$. The excitation

$$z(t;\theta) = \kappa_1 \cdot y(t;\theta) + \kappa_3 \cdot y^3(t;\theta), \tag{2}$$

is a linear-plus-cubic-Gaussian process, having a bimodal first-order pdf in the case $\kappa_1 \cdot \kappa_3 < 0$. Without any loss of generality, we shall assume that $\kappa_1 > 0$. The signs of $\mu_1, \mu_3, \kappa_3$ may be either $+1$, or $-1$, affecting the structure of the solution to both the deterministic and the random problem.

Furthermore, the Gaussian process $y(t;\theta)$ is assumed to be m.s.-continuous, with continuous path functions, and asymptotically stationary. As a consequence the functions $m_y(t)$ and $C_{yy}(t,s)$ are well-defined and bounded in the infinite intervals $[t_0, \infty)$ and $[t_0, \infty) \times [t_0, \infty)$, respectively, and they tend to well-defined limits as $t \to \infty$ and $s \to \infty$. The condition of asymptotic stationarity means that there exist a constant $m_y^{(\infty)}$ and a bounded stationary covariance function $C_{yy}^{(\infty)}(\tau)$, such that

$$(\forall \varepsilon > 0)(\exists t_\varepsilon > t_0): \ [t \geq t_\varepsilon] \Rightarrow \left| m_y(t) - m_y^{(\infty)} \right| < \varepsilon, \tag{3a}$$

and

$$(\forall \varepsilon > 0)(\exists t_\varepsilon > t_0): \ [t \geq t_\varepsilon \wedge s \geq t_\varepsilon] \Rightarrow \left| C_{yy}(t,s) - C_{yy}^{(\infty)}(t-s) \right| < \varepsilon. \tag{3b}$$

Our main goal in this paper is to exploit $m_y(t)$ and $C_{yy}(t,s)$ in order to formulate and solve equations for the response mean value $m_x(t)$ and the two-time covariances $C_{xy}(t,s)$ and $C_{xx}(t,s)$, for $t, s \in [t_0, \infty) \times [t_0, \infty)$. To simplify the data description, we shall assume that $m_y(t) = m_y^{(\infty)}$ for all $t \geq t_0$, and $C_{yy}(t,s) = C_{yy}^{(\infty)}(t-s)$ for all $t \geq t_0$, $s \geq t_0$. Of specific interest herewith is the correlation time $\tau_{yy}^{corr}$ of the excitation processes, which measures the "color" of the processes and quantify their difference from the delta-correlated ones (having $\tau_{yy}^{corr} = 0$). Among the various definitions appearing in the literature we shall follow the one proposed by Stratonovich [61] (see also [62]),

$$\tau_{yy}^{corr} = \frac{1}{\sigma_y^2} \int_0^\infty |C_{yy}^{(\infty)}(u)| du.$$

Numerical results will be presented for the following input processes $y(t;\theta)$:

I. Centered **Ornstein – Uhlenbeck** process (**OU**),
with $m_y^{(\infty)} = 0$, covariance $C_{yy}^{(\infty)}(u) = \sigma_y^2 \cdot e^{-a|u|}$, $a > 0$,

spectral density $S_{yy}^{(\infty)}(\omega) = \frac{\sigma_y^2}{\pi} \frac{a}{a^2 + \omega^2}$, and correlation time $\tau_{yy}^{corr} = 1/a$.

II. Shifted **Ornstein – Uhlenbeck** process, with
$m_y^{(\infty)} = 0$,

$$C_{yy}^{(\infty)}(u) = \sigma_y^2 \cdot e^{-a|u|} \cos(\omega_0 u), \quad a > 0, \tag{4a}$$

$$S_{yy}^{(\infty)}(\omega) = \frac{\sigma_y^2}{2\pi} \left\{ \frac{a}{a^2 + (\omega_0 + \omega)^2} + \frac{a}{a^2 + (\omega_0 - \omega)^2} \right\}, \tag{4b}$$

$$\tau_{yy}^{corr} = \frac{a}{a^2 + \omega_0^2} + \frac{e^{-a\pi/(2\omega_0)}}{1 - e^{-a\pi/(\omega_0)}} \cdot \frac{2\omega_0}{a^2 + \omega_0^2}, \quad \omega_0 > 0. \tag{4c}$$

III. **Gaussian filter** (**Gf**)
with $m_y^{(\infty)} = 0$, $C_{yy}^{(\infty)}(u) = \sigma_y^2 \cdot e^{-au^2}$, $a > 0$,

$$S_{yy}^{(\infty)}(\omega) = \frac{\sigma_y^2}{2\sqrt{\pi a}} \cdot e^{-\omega^2/(4a)}, \quad \tau_{yy}^{corr} = \sqrt{\pi}/(2\sqrt{a}).$$

IV. Shifted **Gaussian filter** with
$m_y^{(\infty)} = 0$,

$$C_{yy}^{(\infty)}(u) = \sigma_y^2 \cdot e^{-au^2} \cos(\omega_0 u), \quad a > 0, \tag{5a}$$

$$S_{yy}^{(\infty)}(\omega) = \frac{\sigma_y^2}{4\sqrt{\pi a}} \cdot \left\{ e^{-(\omega - \omega_0)^2/(4a)} + e^{-(\omega + \omega_0)^2/(4a)} \right\}. \tag{5b}$$

(In the latter case the function $\tau_{yy}^{corr} = \tau_{yy}^{corr}(a, \omega_0)$ cannot be expressed in closed form. It is calculated numerically when needed, as e.g., in Sec. 6). The auto-covariances (4a), (5a), and the spectral density functions (4b), (5b) are plotted in Figs. 1 and 2.

An interesting difference between the OU process (case I) and Gf type processes (cases III and IV) is that the former can be obtained as a solution of a standard Itô SDE, while the later cannot (to the best of the knowledge of the present authors). The latter assertion is further supported by the pertinent discussion at the end of Sec. 5.

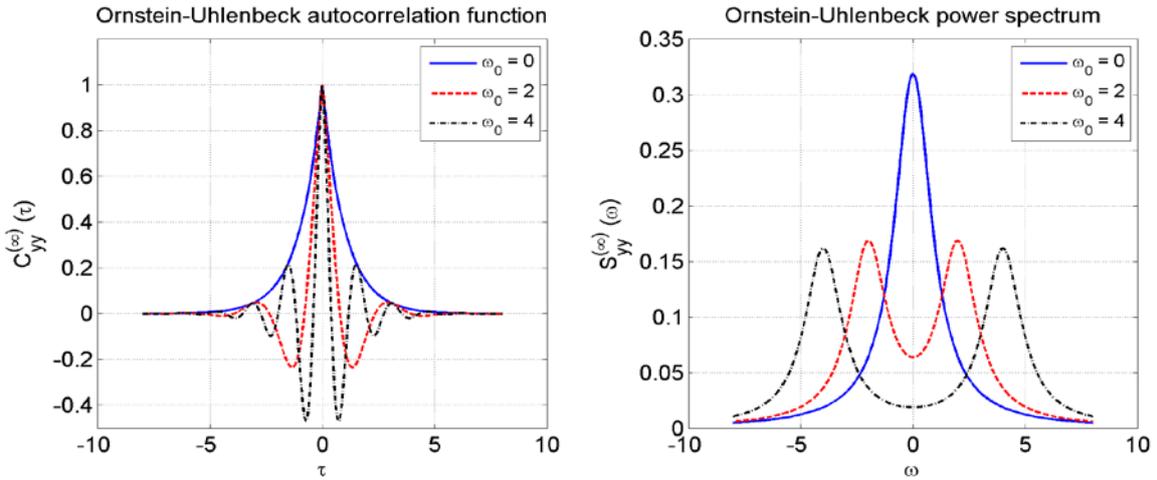

**Fig. 1:** The auto-covariance and the spectral density functions of the considered Ornstein – Uhlenbeck processes

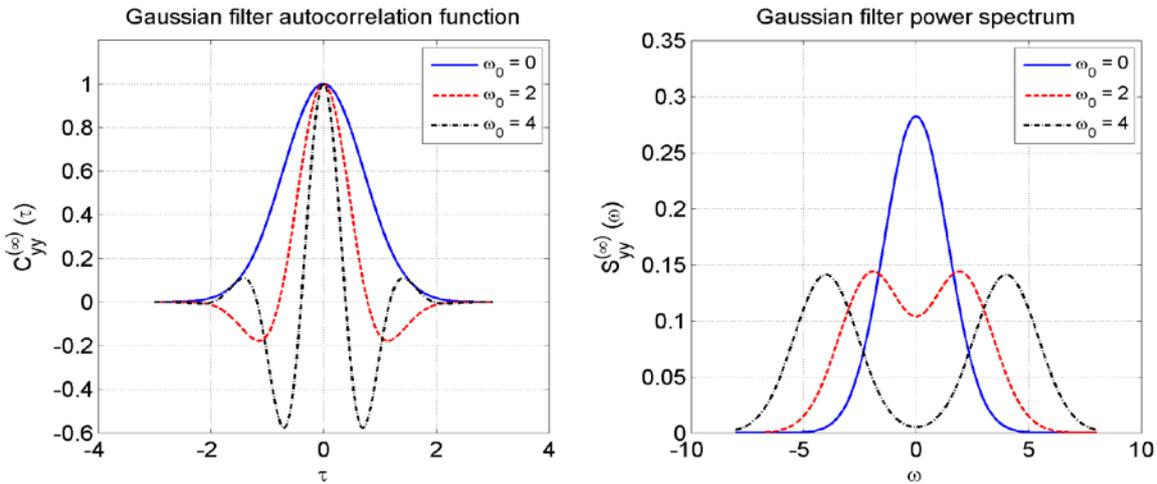

**Fig. 2:** The auto-covariance and the spectral density functions of the considered Gaussian filter processes

Equation (1a) can be rewritten in the form

$$\dot{x}(t) = -\frac{dU(x(t))}{dx} + \kappa_1 \cdot y(t) + \kappa_3 \cdot y^3(t), \tag{6}$$

where $U(x) = -\frac{1}{2}\mu_1 x^2 - \frac{1}{4}\mu_3 x^4$, is a potential function. The shape of $U(x)$ is essentially dependent on the signs of $\mu_1$, $\mu_3$. The four possible configurations of $U(x)$ are shown in Fig. 3. The case $\mu_1 < 0$, $\mu_3 < 0$ (**Fig. 3a**) is referred to as the monostable case and the case $\mu_1 > 0$, $\mu_3 < 0$ (**Fig. 3b**) as the bistable case, in accordance with the stability properties of the homogeneous equation. In both cases equation (6) has the BIBO (bounded-input, bounded-output) property, and thus the long-time limits of solutions under bounded excitations are well defined. In the present work we focus on the (simpler) monostable case. The bistable case will be considered in a forthcoming paper. The cases $\mu_1 < 0$, $\mu_3 > 0$ (**Fig. 3c**) and $\mu_1 > 0$, $\mu_3 > 0$ (**Fig. 3d**) do not exhibit the BIBO property, and thus, in general, they have not finite long-time limits.

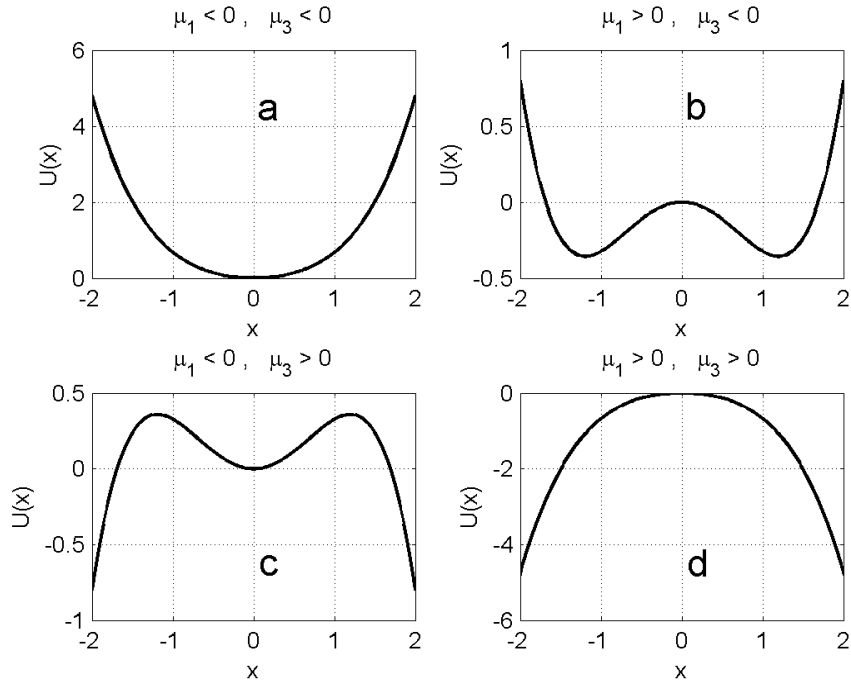

**Fig. 3:** The shape of the potential function $U(x) = -\frac{1}{2}\mu_1 x^2 - \frac{1}{4}\mu_3 x^4$, and the stability characterizations of the homogeneous deterministic problem. **a:** The monostable case, **b:** The bistable case, **c:** Locally stable (at $x = 0$) case, **d:** Globally unstable case

## 3. Derivation of the two-time moment equations

Since the excitation in the RDE (1a,b) is not a delta-correlated process, the correlation structure between $x(t;\theta)$ and $y(s;\theta)$, for $t \neq s \geq t_0$, is important, and should be considered for all $t, s \geq t_0$. We shall now derive and discuss a set of differential equations for the first- and second-order RE moments $m_x(t)$, $C_{xy}(t,s)$, $C_{xx}(t,s)$. By applying the mean value operator to Eq. (1a), we obtain

$$\frac{dm_x(t)}{dt} = \mu_1 \cdot m_x(t) + \mu_3 \cdot R_{xx}^{21}(t,t) + \kappa_1 \cdot m_y(t) + \kappa_3 \cdot R_{yy}^{21}(t,t), \tag{7a}$$

$$m_x(t_0) = m_0. \tag{7b}$$

Multiplying Eqs. (1a,b) first by $y(s;\theta)$, and then by $x(s;\theta)$, and applying the mean value operator, we get the following differential equations and initial conditions for the correlation functions $R_{xy}(t,s)$ and $R_{xx}(t,s)$:

$$\frac{\partial R_{xy}(t,s)}{\partial t} = \mu_1 \cdot R_{xy}(t,s) + \mu_3 \cdot R_{xy}^{31}(t,s) + \kappa_1 \cdot R_{yy}(t,s) + \kappa_3 \cdot R_{yy}^{31}(t,s), \tag{8a}$$

$$R_{xy}(t_0,s) = E^\theta \left[ x_0(\theta) \cdot y(s;\theta) \right] = m_{x_0} \cdot m_y(s), \tag{8b}$$

and

$$\frac{\partial R_{xx}(t,s)}{\partial t} = \mu_1 \cdot R_{xx}(t,s) + \mu_3 \cdot R_{xx}^{31}(t,s) + \kappa_1 \cdot R_{xy}(s,t) + \kappa_3 \cdot R_{xy}^{13}(s,t), \tag{9a}$$

$$R_{xx}(t_0,s) \equiv R_{x_0 x}(s) = E^\theta \left[ x_0(\theta) \cdot x(s;\theta) \right]. \tag{9b}$$

The initial condition (9b) is not known, since it depends on the unknown response $x(s;\theta)$. We have, thus, to derive an equation permitting us to calculate the one-time moment $R_{xx}(t_0,t) = R_{x_0 x}(t)$. Such an equation is easily obtained by multiplying Eq. (1a) by $x_0(\theta)$ and taking mean values:

$$\frac{dR_{x_0 x}(t)}{dt} = \mu_1 \cdot R_{x_0 x}(t) + \mu_3 \cdot R_{x_0 x}^{3}(t) + \kappa_1 \cdot m_y(t) m_{x_0} + \kappa_3 \cdot R_{yy}^{21}(t,t) m_{x_0}. \tag{9c}$$

The initial condition for the latter equation is the known quantity

$$R_{x_0 x}(t_0) = E^\theta \left[ x_0(\theta) \cdot x_0(\theta) \right] = R_{x_0 x_0}. \tag{9d}$$

Since $\partial R_{xx}(t,s)/\partial t \big|_{s=t} \neq dR_{xx}(t,t)/dt$, differential equations (8a) and (9a) cannot be applied (as they stand) to the time-diagonal case $s = t$. It is possible, however, to obtain a differential equation for $R_{xx}(t,t)$. This can be done by multiplying Eq. (1a) by $2x(t;\theta)$ and then applying the mean value operator. The resulting equation and the corresponding initial conditions are

$$\frac{dR_{xx}(t,t)}{dt} = 2 \cdot \mu_1 \cdot R_{xx}(t,t) + 2 \cdot \mu_3 \cdot R_{xx}^{31}(t,t) + 2 \cdot \kappa_1 \cdot R_{xy}(t,t) + 2 \cdot \kappa_3 \cdot R_{xy}^{13}(t,t), \tag{10a}$$

$$R_{xx}(t_0,t_0) = R_{x_0 x_0}. \tag{10b}$$

In order to facilitate the application of the Gaussian moment closure to Eqs. (7) – (10), it is expedient to reformulate them in terms of the central moment functions $C_{xy}^{j_1 j_2}(t,s)$. Relations between $R_{xy}^{j_1 j_2}(t,s)$ and $C_{xy}^{j_1 j_2}(t,s)$ are easily obtained by means of the binomial expansion:

$$R_{xy}^{j_1 j_2}(t,s) = \sum_{k_1=0}^{j_1} \sum_{k_2=0}^{j_2} \binom{j_1}{k_1} \binom{j_2}{k_2} m_x^{j_1-k_1}(t) \, m_y^{j_2-k_2}(s) \, C_{xy}^{k_1 k_2}(t,s). \tag{11}$$

Using appropriate special cases of Eq. (11), Eqs. (7) – (10) are transformed to the following equivalent forms:

$$\frac{dm_x(t)}{dt} = \mu_1 \cdot m_x(t) + \mu_3 \cdot m_x^3(t) + \mu_3 \cdot \left(C_{xx}^{21}(t,t) + 3m_x(t) C_{xx}(t,t)\right) +$$
$$+ \kappa_1 \cdot m_y(t) + \kappa_3 \cdot m_y^3(t) + \kappa_3 \cdot \left(C_{yy}^{21}(t,t) + 3m_y(t) C_{yy}(t,t)\right), \tag{12a}$$

$$m_x(t_0) = m_{x_0}. \tag{12b}$$

$$\frac{\partial C_{xy}(t,s)}{\partial t} = \mu_1 \cdot C_{xy}(t,s) + \kappa_1 \cdot C_{yy}(t,s) +$$
$$+ \mu_3 \cdot \left(C_{xy}^{31}(t,s) + 3m_x(t) C_{xy}^{21}(t,s) + 3m_x^2(t) C_{xy}(t,s)\right) + \tag{13a}$$
$$+ \kappa_3 \cdot \left(C_{yy}^{31}(t,s) + 3m_y(t) C_{yy}^{21}(t,s) + 3m_y^2(t) C_{yy}(t,s)\right),$$

$$C_{xy}(t_0,s) = 0. \tag{13b}$$

$$\frac{\partial C_{xx}(t,s)}{\partial t} = \mu_1 \cdot C_{xx}(t,s) + \kappa_1 \cdot C_{xy}(s,t) +$$
$$+ \mu_3 \cdot \left(C_{xx}^{31}(t,s) + 3m_x(t) C_{xx}^{21}(t,s) + 3m_x^2(t) C_{xx}(t,s)\right) + \tag{14a}$$
$$+ \kappa_3 \cdot \left(C_{xy}^{13}(s,t) + 3m_y(t) C_{xy}^{12}(s,t) + 3m_y^2(t) C_{xy}(s,t)\right),$$

$$C_{xx}(t_0,s) \equiv C_{x_0 x}(s) = E^\theta \left[x_0(\theta) \cdot x(s;\theta)\right] - m_{x_0} \cdot m_x(s), \tag{14b}$$

$$\frac{dC_{x_0 x}(t)}{dt} = \left(\mu_1 + 3\mu_3 m_x^2(t) \cdot\right) C_{x_0 x}(t) + \mu_3 \cdot \left(C_{x_0 x}^{3}(t) + 3m_x(t) C_{x_0 x}^{2}(t)\right), \tag{14c}$$

$$C_{x_0 x}(t_0) = C_{x_0 x_0}, \tag{14d}$$

$$\frac{dC_{xx}(t,t)}{dt} = 2\cdot\left(\mu_1 + 3\mu_3 m_x^2(t)\right)\cdot C_{xx}(t,t) + 2\cdot\left(\kappa_1 + 3\kappa_3 m_y^2(t)\right)\cdot C_{xy}(t,t) + \tag{15a}$$
$$+ 2\cdot\mu_3\cdot\left(C_{xx}^{31}(t,t) + 3m_x(t) C_{xx}^{21}(t,t)\right) + 2\cdot\kappa_3\cdot\left(C_{xy}^{13}(t,t) + 3m_y(t) C_{xy}^{12}(t,t)\right),$$
$$C_{xx}(t_0,t_0) = C_{x_0 x_0}. \tag{15b}$$

The five differential equations (12a), (13a), (14a), (14c) and (15a) will be considered as a system of equations for the five moment functions $m_x(t)$, $C_{xy}(t,s)$, $C_{xx}(t,s)$, $C_{xx}(t,t)$, $C_{x_0 x}(t)$. Clearly, these equations are not closed in two ways. On the one hand, they contain higher-order moments, like $C_{xx}^{31}(t,s)$, $C_{xx}^{21}(t,s)$, etc., as also the usual one-time moments do. On the other hand, two of them, namely Eqs. (13a) and (14a), being differential equations only with respect to the response time $t$, contain a second time variable $s$, acting as a parameter. This fact renders them not closed in time as well. In the next section we shall elaborate closure schemes in both respects, obtaining a closed and solvable system of differential equations, non-linear and non-local in time.

Before closing this section we shall comment on an apparently controversial feature of Eq. (14a), describing the evolution of $C_{xx}(t,s)$. In this equation, terms containing $C_{xy}(s,t)$, $C_{xy}^{12}(s,t)$ and $C_{xy}^{13}(s,t)$ are included in the last term of right-hand side, bringing into play an apparently non-causal dependence: the effect of $y(t;\theta)$ on $x(s;\theta)$, for $s<t$. This apparent controversy is resolved in the next section, by appropriately combining Eqs. (14a) and (13a), resulting in a dependence of $C_{xx}(t,s)$ upon the whole history of the data moment $C_{yy}(\tau,t)$, for all $\tau\in[t_0,t]$. Thus, there is not a violation of causality, but simply a manifestation of the non-locality in time of the above equations, due to the smooth time dependence of $C_{yy}(t,s)$, $C_{xy}(t,s)$ and $C_{xx}(t,s)$ for all $t,s\geq t_0$. As a result, all three functions "live" (are supported) in the whole rectangle

$ABCD: R_{ts}(T) = \left\{(t,s): t_0 \leq t \leq T,\ t_0 \leq s \leq T\right\}$,

not only in the triangle $ABC$ (**Fig. 4**).

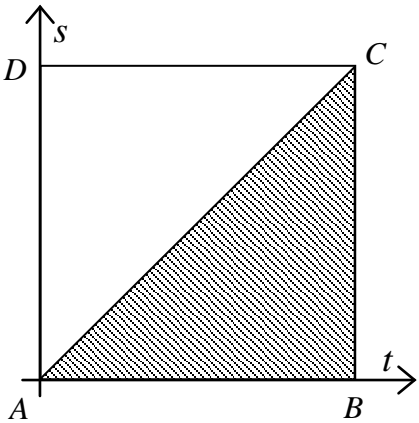

**Fig. 4.** The two-time rectangle $R_{ts}$

## 4. Closure schemes of the two-time moment Equations

### 4.1 Moment closure

To eliminate higher-order moments from the right-hand side of equations (12a), (13a), (14a), (14c) and (15a), use will be made of the Gaussian closure assumption. Under this assumption all the third-order central moments vanishes, and the forth-order ones are expressed by means of second-order central moments, in accordance to Isserlis' Theorem:

$$C_{xx}^{31}(t,s) = 3 \cdot C_{xx}(t,t) \cdot C_{xx}(t,s), \tag{16a}$$

$$C_{xy}^{31}(t,s) = 3 \cdot C_{xx}(t,t) \cdot C_{xy}(t,s), \tag{16b}$$

$$C_{xy}^{13}(s,t) = 3 \cdot C_{yy}(t,t) \cdot C_{xy}(s,t), \tag{16c}$$

$$C_{yy}^{31}(t,s) = 3 \cdot C_{yy}(t,t) \cdot C_{yy}(t,s), \tag{16d}$$

$$C_{x_0 x}^{3}(t) = 3 \cdot C_{xx}(t,t) \cdot C_{x_0 x}(t). \tag{16e}$$

Introducing the approximations (16) into Eqs. (12a), (13a), (14a), (14c) and (15a), we obtain:

$$\frac{dm_x(t)}{dt} = \left(\mu_1 + \mu_3 m_x^2(t) + 3\mu_3 C_{xx}(t,t)\right) \cdot m_x(t) + \tilde{B}_y(t) \cdot m_y(t), \tag{17}$$

$$\frac{\partial C_{xy}(t,s)}{\partial t} = \left(\mu_1 + 3\mu_3 m_x^2(t) + 3\mu_3 C_{xx}(t,t)\right) C_{xy}(t,s) + B_y(t) \cdot C_{yy}(t,s), \tag{18}$$

$$\frac{\partial C_{xx}(t,s)}{\partial t} = \left(\mu_1 + 3\mu_3 m_x^2(t) + 3\mu_3 C_{xx}(t,t)\right) C_{xx}(t,s) + B_y(t) \cdot C_{xy}(s,t), \tag{19}$$

$$\frac{dC_{x_0 x}(t)}{dt} = \left(\mu_1 + 3\mu_3 m_x^2(t) + 3\mu_3 C_{xx}(t,t)\right) C_{x_0 x}(t), \tag{20}$$

$$\frac{dC_{xx}(t,t)}{dt} = 2 \cdot \left(\mu_1 + 3\mu_3 m_x^2(t) + 3\mu_3 C_{xx}(t,t)\right) \cdot C_{xx}(t,t) + 2 \cdot B_y(t) \cdot C_{xy}(t,t), \tag{21}$$

where

$$\tilde{B}_y(t) = \kappa_1 + \kappa_3 m_y^2(t) + 3\kappa_3 C_{yy}(t,t), \tag{22a}$$

$$B_y(t) = \kappa_1 + 3\kappa_3 m_y^2(t) + 3\kappa_3 C_{yy}(t,t), \tag{22b}$$

Having been working off the higher-order moments, it remains to elaborate on the second special feature of the two-time moment equations, the simultaneous appearance of the two time variables $t$, $s$.

### 4.2 Time closure

We shall now implement a time closure, obtaining an one-time, closed, causal subsystem of two equations for the time-diagonal moments $m_x(t)$ and $C_{xx}(t,t)$. Setting

$$A_x(t) \equiv A_x\left[m_x(t), C_{xx}(t,t)\right] = \mu_1 + 3\mu_3 m_x^2(t) + 3\mu_3 C_{xx}(t,t), \tag{22c}$$

we rewrite equation (18) in the form

$$\frac{\partial C_{xy}(t,s)}{\partial t} = A_x(t) \cdot C_{xy}(t,s) + B_y(t) \cdot C_{yy}(t,s). \tag{23}$$

Although the function $A_x(t)$ is not known, it possesses two important properties: First, it is not dependent on the $C_{xy}(t,s)$ itself and, second, in the monostable case studied herewith ($\mu_1, \mu_3 < 0$), it is always negative:

$$A_x(t) < 0, \quad \text{for all } t \geq t_0. \tag{24}$$

On the basis of the first property, we can consider Eq. (23) as a linear, first-order ODE for the cross-covariance $C_{xy}(t,s)$ with respect to $t$ ($s$ being considered as a parameter). In accordance with the standard theory of first-order ODEs [63], the solution of Eq. (23) with initial condition (13b), can be expressed by the following integral formula, in terms of the unknown function $A_x(t)$:

$$C_{xy}(t,s) = \int_{t_0}^{t} B_y(\tau) \, C_{yy}(\tau,s) \exp\left(\int_{\tau}^{t} A_x(u) \, du\right) d\tau. \tag{25}$$

This solution is valid for any $t \geq t_0$ and for any $s \neq t$. However, looking closer at the structure of the right-hand side of Eq. (25), we observe that it is a continuous function on $s$ for all $s \geq t_0$ (since it depends on $s$ only through the continuous data function $C_{yy}(\tau,s)$). Thus, taking the limit of both sides of Eq. (25) for $s \to t$, we obtain

$$C_{xy}(t,t) = \int_{t_0}^{t} B_y(\tau) \, C_{yy}(\tau,t) \exp\left(\int_{\tau}^{t} A_x\left[m_x(u), C_{xx}(u,u)\right] du\right) d\tau. \tag{26a}$$

Eq. (26a) expresses the time-diagonal cross-covariance $C_{xy}(t,t)$ as a non-linear, causal operator on the whole history of the (unknown) response mean value $m_x(u)$, $t_0 \leq u \leq t$, and response, time-diagonal auto-covariance $C_{xx}(u,u)$, $t_0 \leq u \leq t$. This is the key point for the time closure.

**Remark:** By direct differentiation it is easy to verify that $C_{xy}(t,s)$, as expressed by Eq. (25), satisfies Eq. (23). What is more interesting here is to use (26a) in order to find a differential equation for the time-diagonal cross-covariance $C_{xy}(t,t)$. Again by direct differentiation, of Eq. (26a) this time, we find

$$\frac{dC_{xy}(t,t)}{dt} = A_x(t) \cdot C_{xy}(t,t) + B_y(t) \cdot C_{yy}(t,t) +$$

$$+ \int_{t_0}^{t} B_y(\tau) \frac{\partial C_{yy}(\tau,t)}{\partial t} \exp\left(\int_{\tau}^{t} A_x(u)\, du\right) d\tau. \quad (26b)$$

The latter is a non-linear and non-local in time (causal) equation. (Recall that $A_x(t)$ is an unknown function, implicitly dependent on $C_{xy}(t,t)$, because of Eq. (21)). In Sec. 5 we shall see how this non-local equation is localized, if the data random function $y(t;\theta)$ is an OU process which is derived by means of an Itô SDE.

Substituting $C_{xy}(t,t)$ from Eq. (26a) in Eq. (21), we obtain

$$\frac{dC_{xx}(t,t)}{dt} = 2 \cdot \left(\mu_1 + 3\mu_3\, m_x^2(t) + 3\mu_3\, C_{xx}(t,t)\right) \cdot C_{xx}(t,t) +$$

$$+ 2 \cdot B_y(t) \cdot \int_{t_0}^{t} B_y(\tau)\, C_{yy}(\tau,t) \exp\left(\int_{\tau}^{t} \left(\mu_1 + 3\mu_3\, m_x^2(u) + 3\mu_3\, C_{xx}(u,u)\right) du\right) d\tau. \quad (27)$$

Eqs. (17) and (27), with initial conditions (12b) and (15b), respectively, form a closed, non-linear, causal system of evolution equations for the moment functions $m_x(t), C_{xx}(t,t)$.

### 4.3  *On the solvability of the causal moment system* (17), (27)

The system (17), (27) belongs to the family of causal differential systems, for which an extensive literature has been developed in the last two decades; see the books [46], [47], [48], [49] and the references cited there. Local existence and uniqueness for system (17), (27) can be easily established by invoking Theorem 2.3.2 of [47]. Global existence in $[t_0, +\infty)$ can be proved as well, on the basis of Theorem 2.4.1 of the same book, but the treatment is much more complicated. The existence of limiting values $m_x^{(\infty)} = \lim_{t \to \infty} m_x(t)$ and $C_{xx}^{(\infty)}(0) = \lim_{t \to \infty} C_{xx}(t,t)$ is ensured on the basis of the property (24). A robust and efficient, iterative, numerical solution to the causal system (17), (27), in the whole time domain $[t_0, +\infty)$, is presented in Sec. 6.

### 4.4  *Representation of the off-diagonal moments $C_{xy}(t,s)$ and $C_{xx}(t,s)$ in terms of the diagonal ones*

Having solved the system (17), (27) and obtained the functions $m_x(t)$ and $C_{xx}(t,t)$, Eq. (25) provides us with the two-time cross-covariance $C_{xy}(t,s)$ in the whole rectangle $R_{ts}$ (**Fig. 4**). Furthermore, substituting $C_{xy}(s,t)$ in the right-hand side of Eq. (19) from (25) (with arguments $t$ and $s$ interchanged), Eq. (19) becomes a first-order ODE with known variable coefficients of the form

$$\frac{\partial C_{xx}(t,s)}{\partial t} = A_x(t)\,C_{xx}(t,s) + F_{xy}(t,s), \tag{28a}$$

where

$$F_{xy}(t,s) = \int_{\tau=t_0}^{\tau=s} B_y(\tau)\,C_{yy}(\tau,t)\,B_y(t)\exp\left(\int_{u=\tau}^{u=s} A_x(u)\,du\right) d\tau. \tag{28b}$$

For any $s \geq t_0$, the initial condition $C_{xx}(t_0,s) = C_{x_0 x}(s)$ is calculated by solving Eq. (20) with initial condition $C_{xx}(t_0,t_0) = C_{x_0 x_0}$:

$$C_{x_0 x}(s) = C_{x_0 x_0} \cdot \exp\left(\int_{t_0}^{s} A_x(u)\,du\right). \tag{28c}$$

The solution of Eq. (28a) with initial condition (28c) is given by the formula

$$C_{xx}(t,s) = C_{x_0 x_0} \cdot \exp\left(\int_{t_0}^{t} A_x(u)\,du + \int_{t_0}^{s} A_x(u)\,du\right) +$$

$$+ \int_{\tau_1=t_0}^{\tau_1=t} \int_{\tau_2=t_0}^{\tau_2=s} G_{yy}(\tau_1,\tau_2) \exp\left(\int_{u_1=\tau_1}^{u_1=t} A_x(u_1)\,du_1\right) \exp\left(\int_{u_2=\tau_2}^{u_2=s} A_x(u_2)\,du_2\right) d\tau_1 d\tau_2,$$
(29a)

where

$$G_{yy}(\tau_1,\tau_2) = B_y(\tau_1)\,C_{yy}(\tau_1,\tau_2)\,B_y(\tau_2). \tag{29b}$$

The symmetry relation $C_{xx}(t,s) = C_{xx}(s,t)$ is clearly revealed by Eqs. (29).

## 5. Equivalence of causal moment equations with those obtained by the Itô approach for a centered OU stochastic input

It is well-known that the OU process $y(t;\theta)$, with $m_y(t) = m_y^{(\infty)} = 0$ and $C_{yy}(t,s) = C_{yy}^{(\infty)}(t-s) = \sigma_y^2 \cdot \exp(-a|t-s|)$ is the solution to the Itô SDE

$$\dot{y}(t;\theta) = -a\,y(t;\theta) + \xi(t;\theta), \tag{30}$$

where $\xi(t;\theta)$ is a Gaussian with noise, with $\sigma_\xi^2 = 2a\sigma_y^2$. Thus, in this case, it is possible to consider Eq. (30), along with (1a), as a system of two Itô SDEs, and derive the moment equations from this. The corresponding equations for the (one-time) first- and second-order moments can be obtained by using either the FPK equation or the Itô Lemma, and they have the form

$$\frac{dm_x(t)}{dt} = \mu_1 \cdot m_x(t) + \mu_3 \cdot R_{xx}^{21}(t,t) + \kappa_3 \cdot R_{yy}^{21}(t,t), \tag{31a}$$

$$\frac{dR_{xy}(t,t)}{dt} = (\mu_1 - a) \cdot R_{xy}(t,t) + \mu_3 \cdot R_{xy}^{31}(t,t) + \kappa_1 \cdot R_{yy}(t,t) + \kappa_3 \cdot R_{yy}^{31}(t,t), \tag{31b}$$

$$\frac{dR_{xx}(t,t)}{dt} = 2\mu_1 \cdot R_{xx}(t,t) + 2\mu_3 \cdot R_{xx}^{31}(t,t) + 2\kappa_1 \cdot R_{xy}(t) + 2\kappa_3 \cdot R_{xy}^{13}(t,t). \tag{31c}$$

After transforming Eqs. (31a,b,c) to equivalent ones with respect to the central moments, and applying Gaussian closure, we obtain

$$\frac{dm_x(t)}{dt} = \left(\mu_1 + \mu_3 m_x^2(t) + 3\mu_3 C_{xx}(t,t)\right) \cdot m_x(t), \tag{32a}$$

$$\frac{dC_{xy}(t,t)}{dt} = \left(\mu_1 - a + 3\mu_3 C_{xx}(t,t) + 3\mu_3 m_x^2(t)\right) \cdot C_{xy}(t,t) + B_y(t) \cdot C_{yy}(t,t) \tag{32b}$$

$$\frac{dC_{xx}(t,t)}{dt} = 2 \cdot \left(\mu_1 + 3\mu_3 m_x^2(t) + 3\mu_3 C_{xx}(t,t)\right) \cdot C_{xx}(t,t) + 2 \cdot B_y(t) \cdot C_{xy}(t,t). \tag{32c}$$

This is a closed system of three non-linear ODEs, which can be solved numerically. The most interesting question arising at this point is: *what is the relation of the closed system* (32) *with the moment equations derived in the previous section*. Eqs. (32a) and (32c) are identical with our Eqs. (17) (with $m_y(t) = 0$) and (21). Eq. (32b) is different from our Eq. (18), as expected, since the latter refers to the two-time covariance $C_{xy}(t,s)$. It turns out that Eq. (32b) should be compared with our Eq. (26b), which also refers to $C_{xy}(t,t)$. For the specific choice considered, where $C_{yy}(t,s) = \sigma_y^2 \cdot \exp(-a|t-s|)$, we have

$$\left. \frac{\partial C_{yy}(\tau,t)}{\partial t} \right|_{\tau \leq t} = -a\, C_{yy}(\tau,t). \tag{33}$$

Substituting the derivative $\partial C_{yy}(\tau,t)/\partial t$ in (the non-local term of) Eq. (26b) from Eq. (33), we find that Eq. (26b) becomes identical with Eq. (32b). This means that, in the specific case of a centered OU process $y(t;\theta)$, for which Eq. (33) holds true, system (32) is equivalent with our causal system (17), (27) and (26b). This equivalence shows that our approach consistently generalizes the Itô/filtering approach, remaining valid for any kind of covariance $C_{yy}(t,s)$. It is also instructive to observe that Eq. (26b) is not localized if $y(t;\theta)$ is a shifted OU or a Gf process. This can be considered as an indication that the latter processes are not representable as solutions of Itô SDEs. This assertion is posed here simply as a conjecture, since the present authors are not aware of any rigorous discussion or proof of the issue.

## 6. Numerical solution to the two-time moment equations

### 6.1 Method of numerical solution

The numerical solution to the two-time moment system (17) – (21) is performed in two stages. In the first stage the one-time, non-local (causal) system of Eqs. (17), (27) is solved, providing the moments $m_x(t)$ and $C_{xx}(t,t)$, from which the function $A_x(t)$ is calculated as well (using Eq. (22c)). The knowledge of $A_x(t)$ allows us to proceed to the second stage, that is, the numerical calculation of the two-time moments $C_{xy}(t,s)$ and $C_{xx}(t,s)$. This can be done either by using the integral forms (25) and (29a), or by integrating the corresponding ODEs (23) and (28a).

- *First stage: Numerical solution to the one-time causal system* **(17), (27)**

Since Eq. (27) contains a causal term dependent on the history of both unknowns $m_x(t)$ and $C_{xx}(t,t)$ over all past times, special techniques should be invoked in order to solve numerically the system (17), (27). The numerical method used in this work is a two-scale, iterative prediction/correction scheme, alternating between the prediction of the function $A_x(t)$ (needed to calculate the non-local term in the coarse scale), and the solution of a usual nonlinear system of ODEs, which provides us with new, corrected versions of $m_x(t)$, $C_{xx}(t,t)$ and $A_x(t)$, as well (in the fine scale). The numerical solution to ODEs is performed by means of the adaptive fourth order Runge-Kutta method developed by Dormant and Prince [64], implemented in MATLAB by the ODE45 function.

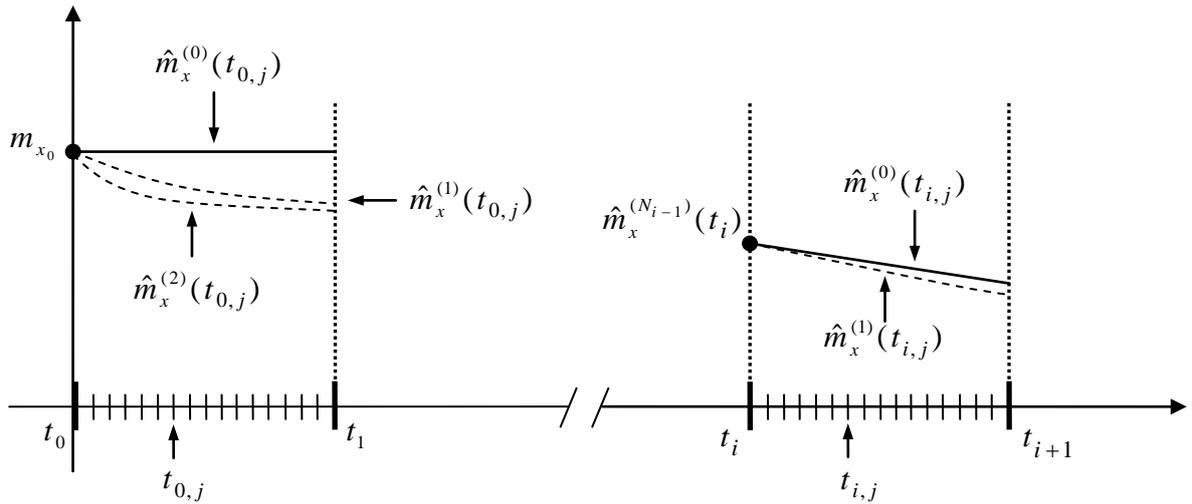

**Fig. 5.** Coarse-scale ($t_0, t_1, \ldots, t_i, \ldots$) and fine-scale ($t_{i,0}, t_{i,1}, \ldots, t_{i,j}, \ldots$) discretization. Various cycles approximations are depicted for the mean value $m_x(t)$. The covariance $C_{xx}(t,t)$ is approximated in a similar way.

To explain better the method of numerical solution, we first introduce some specific notation and terminology. Let $[t_0, T] = [0, T]$ be a time interval, long enough so that the long-time statistical equilibrium limit has been reached. This interval is discretized by a (not necessarily uniform) discretization $\mathcal{D}^{\text{coarse}} = \{t_i : 0 = t_0 < t_1 < ... < t_I = T\}$, called the coarse-scale discretization. Each subinterval $[t_i, t_{i+1}]$ of the coarse-scale discretization is further discretized by intermediate time points $t_{i,j}$, $j = 0(1)J_i$, defining the fine-scale discretization $\mathcal{D}^{\text{fine}}_{(i)} = \{t_{i,j} : t_{i,0} = t_i < t_{i,1} < ... < t_{i,J_i} = t_{i+1}\}$; see **Fig. 5**. The discrete numerical approximations of $m_x(t)$ and $C_{xx}(t,t)$ are denoted by $\hat{m}_x^{(n)}(t_i)$ and $\hat{C}_{xx}^{(n)}(t_i, t_i)$ in the coarse-scale discretization, and by $\hat{m}_x^{(n)}(t_{i,j})$ and $\hat{C}_{xx}^{(n)}(t_{i,j}, t_{i,j})$ in the fine-scale discretization. The upper index $n = 0(1)N_i$ denotes the number of local iterations performed at each time step $[t_i, t_{i+1}]$, to improve the approximation, as described below.

The solution procedure is initialized by using the known initial values $m_x(t_0) = m_{x_0}$ and $C_{xx}(t_0, t_0) = C_{x_0 x_0}$. Assuming that these values remain constant in the interval $[t_0, t_1]$, we formulate the initial approximation of the quantity $A_x(t)$, $A_x^{(0)}(t) = A_x(t_0) = constant$, in $[t_0, t_1]$, using Eq. (22c). Substituting $A_x^{(0)}(t)$ in the non-local term of the right-hand side of Eq. (27), we get the following system of nonlinear ODEs

$$\frac{dm_x^{(0)}(t)}{dt} = \left(\mu_1 + \mu_3 \left(m_x^{(0)}(t)\right)^2 + 3\mu_3 C_{xx}^{(0)}(t,t)\right) \cdot m_x^{(0)}(t) + \tilde{B}_y(t) \cdot m_y(t), \qquad (34a)$$

$$\frac{dC_{xx}^{(0)}(t,t)}{dt} = 2 \cdot \left(\mu_1 + 3\mu_3 \left(m_x^{(0)}(t)\right)^2 + 3\mu_3 C_{xx}^{(0)}(t,t)\right) \cdot C_{xx}^{(0)}(t,t) + \\ + 2 \cdot B_y(t) \cdot \int_{t_0}^{t} B_y(\tau) C_{yy}(\tau, t) \exp\left(\int_{\tau}^{t} A_x^{(0)}(u) \, du\right) d\tau, \qquad (34b)$$

which is solved numerically in $[t_0, t_1]$, in the fine-scale discretization $\mathcal{D}^{\text{fine}}_{(0)}$, generating the zeroth-cycle approximations $\hat{m}_x^{(0)}(t_{0,j})$, $\hat{C}_{xx}^{(0)}(t_{0,j}, t_{0,j})$, $j = 1(1)J_0$. Using these values and Eq. (22c), the function $A_x(t)$ is updated (corrected), obtaining the first-cycle approximation $A_x^{(1)}(t)$ by linear interpolation, $A_x^{(1)}(t) = \sum_{j=0}^{J_0} \hat{A}_x^{(0)}(t_{0,j}) \cdot \ell_j(t)$, where $\ell_j(t)$ are linear interpolants ($\ell_j(t_{0,j'}) = \delta_{j,j'}$, $\text{supp}(\ell_j) = [t_{0,j-1}, t_{0,j+1}]$), and $\hat{A}_x^{(0)}(t_{0,j})$ are obtained by using Eq. (22c). Substituting $A_x^{(0)}(t)$ by $A_x^{(1)}(t)$, in the non-local term of the right-hand side of Eq. (34b), we obtain an updated (corrected) version of the nonlinear system (34), which is again solved numerically, as previously, providing us with the first-cycle approximations $\hat{m}_x^{(1)}(t_{0,j})$, $\hat{C}_{xx}^{(1)}(t_{0,j}, t_{0,j})$, $j = 1(1)J_0$. The procedure is repeated until two successive cycle calculations remain desirably close, i.e.

$$\left| \hat{m}_x^{(N_0-1)}(t_0^{(j)}) - \hat{m}_x^{(N_0)}(t_0^{(j)}) \right| < \varepsilon_1, \qquad j = 1(1)J_0$$
$$\left| \hat{C}_{xx}^{(N_0-1)}(t_0^{(j)}, t_0^{(j)}) - \hat{C}_{xx}^{(N_0)}(t_0^{(j)}, t_0^{(j)}) \right| < \varepsilon_2, \quad j = 1(1)J_0 \tag{35}$$

where $\varepsilon_1$, $\varepsilon_2 > 0$ are user-defined tolerances.

When the aforementioned convergence condition has been satisfied, the numerical scheme can advance in time, to solve the system (17), (27) in $[t_1, t_2]$, providing the approximations $\hat{m}_x^{(N_1)}(t_1^{(j)})$, $\hat{C}_{xx}^{(N_1)}(t_1^{(j)}, t_1^{(j)})$, $j = 0(1)J_1$. In this case, the zeroth-cycle function $A_x^{(0)}(t)$ in $[t_1, t_2]$, need not to be taken as constant. It can be estimated by linearly extrapolating the two last available values, $\hat{A}_x^{(N_0)}(t_0^{(J_0-1)})$ and $\hat{A}_x^{(N_0)}(t_0^{(J_0)}) = \hat{A}_x^{(N_0)}(t_1)$, from the solution of the previous coarse-scale time step. The scheme then continues in the same way until the desired end time $t_I = T$ has been reached.

The one-time cross-covariance $C_{xy}(t, t)$ is calculated numerically by using Eq. (26a). In fact its calculation has essentially been already done as part of the evaluation of the non-local term of Eq. (27).

The numerical behavior of the above describe solution scheme is very robust, converging quickly with respect to the non-local term, usually within two or three cycles in each coarse-scale time step, when using tolerance $\varepsilon_1 = \varepsilon_2 = 10^{-6}$ and discretization step for $\mathcal{D}^{\text{coarse}}$ of the same order of magnitude as the correlation time $\tau_{yy}^{\text{corr}}$ of the input process.

- *Second stage: Calculation of the two-time moments* ($t \neq s$)

Moving on to the calculation of the off-diagonal, second-order moments over the rectangle $R_{ts}(T) = \{(t, s): t_0 \leq t \leq T, t_0 \leq s \leq T\}$, we can work in two, essentially equivalent, ways, both requiring the one-time moment approximations, found in the previous stage. The discrete numerical approximations of $C_{xx}(t, s)$ and $C_{xy}(t, s)$ are denoted by $\hat{C}_{xx}(t, s)$ and $\hat{C}_{xy}(t, s)$.

**One approach** is to use the integral forms, Eqs. (25), (29a), and integrate numerically for each $(t, s)$ pair of the rectangle $R_{ts}(T)$. Note that $\hat{C}_{xx}(t, s)$ suffices to be calculated only on the half part of $R_{ts}(T)$ (e.g. the triangle $ABC$, **Fig. 4**), due to the symmetry of $C_{xx}(t, s)$.

**Another approach**, is to numerically solve the ODEs for the two-time moments, i.e. equations (23), (28a), formulating initial value problems with respect to $t$, for each value of $s$. Each of these problems can be solved starting either form $t = t_0$, with initial conditions $C_{xy}(t_0, s) = 0$ and $C_{xx}(t_0, s) = C_{x_0 x}(s)$ (calculated from Eq. (28c)), or from the diagonal $t = s$, using the previous-stage, one-time solutions $\hat{C}_{xx}^{(n)}(t_{i,j}, t_{i,j})$, $\hat{C}_{xy}^{(n)}(t_{i,j}, t_{i,j})$ as initial values. Both approaches yield practically identical results when using the adaptive Runge-Kutta solver, mentioned earlier.

Both ODE and integral form approaches return nearly identical results. The integral form approach is more robust (with respect to the domain discretization), compared to the ODE approach, which requires cautious numerical tolerance configuration to avoid fluctuating artifacts on the results. It has also the advantage that it can be used to obtain the moments at (any) point of $R_{ts}(T)$, without requiring (as do the ODEs) the solution to start from $t = t_0$ or $t = s$ (when solving from the diagonal). On the other hand, the ODE approach is somewhat more generic, since it does not require the analytical calculation of the integral forms.

### 6.2 Monte Carlo simulation

The full problem of the monostable half-oscillator with cubic damping and Gaussian-plus-cubic-Gaussian excitation, as defined by Eqs. (1a,b), has also been thoroughly investigated using Monte Carlo (MC) simulation. The numerical experiments and statistical analysis of the corresponding data, allow moment approximations beyond the Gaussian closure assumptions and thus, an impeccable assessment of the closure induced error. Furthermore, the simulation serve as a means of verification of the proposed numerical scheme when the deviation from the Gaussian case is small.

The simulation is based on the numerical generation of excitation sample functions and the solution of the differential Eqs. (1a,b). The random generator used for the simulation of the Gaussian process $y(t;\theta)$ was coded in MATLAB and is based on the random-phase model. Descriptions of this model can be found in various books [3], [6], [65]; see also [66], [67].

For the MC simulation of Eq. (1), a code named MCIVP_SIM has been developed (also in MATLAB), which allows the simulation of ODE initial value problem, for any user defined, smoothly-correlated, random excitation and system coefficients, and any random initial condition. The simulator uses the ODE solver set available with MATLAB, each time choosing a solver depending on the stiffness of the system, and also allows the definition of custom solvers. MCIVP_SIM can be executed in parallel on shared memory systems with a nearly linear speed-up with respect to the available processors.

For all results shown in this paper, the initial value was considered to be a random variable with $m_{x_0} = 2$ and $C_{x_0 x_0} = 1$. The simulations were executed using $10^4$ samples for the transient moment results.

## 6.3 Numerical results

Numerical results for the three moment functions $m_x(t)$, $C_{xy}(t,s)$ and $C_{xx}(t,s)$ are presented for a cubic, monostable half-oscillator with $\mu_1 = -1.0$, $\mu_3 = 0.0, -0.1, -0.7$ (linear, mildly nonlinear, strongly nonlinear, respectively), $\kappa_1 = 1.0$, and $\kappa_3 = 0.0, 0.4, -0.4$ (Gaussian excitation, unimodal non-Gaussian excitation, bimodal non-Gaussian excitation, respectively). The relaxation time of this oscillator can be taken as $\tau^{\text{relx}} = 1/|\mu_1| = 1$ (the same as in the linear case $\mu_3 = 0$), since the negative $\mu_3$ value (monostability) just enhances the rate towards the long-time statistical equilibrium. The data random function $y(t;\theta)$, defining the excitation, ranges through various cases of OU and Gf processes, as described in Sec. 2. Recall that, in all such cases $C_{yy}(t,s) = C_{yy}^{(\infty)}(t-s)$ is fully defined by three parameters, the variance $\sigma_y^2$, the correlation time $\tau_{yy}^{\text{corr}}$, and the central spectral frequency $\omega_0$. Of course the shape of $C_{yy}(t,s)$ is different in OU and Gf cases. Transient MC simulations have also been conducted for $10^4$ sample functions, and used for the direct numerical calculation of the evolution of the moment functions, for comparison purposes. In most cases the results are presented in the time interval $[t_0, T] = [0, 3]$, which is long enough in order to reveal the long-time statistical equilibrium limit.

The aim of the results presented in this section is twofold. First, to establish the validity of the new, causal (non-local in time) moment equations, derived in Sec. 4; to investigate their limitations, coming from the use of the Gaussian closure; and to ascertain the correctness of the numerical method and codes developed for their solution. Second, to investigate the effect of the variation of the shape (Gf vs. OU) of the given covariance $C_{yy}(t,s)$ on the response moments, when the defining integral parameters ($\sigma_y^2, \tau_{yy}^{\text{corr}}, \omega_0$) are the same. Among all cases presented herewith, only the one corresponding to a centered OU data function $y(t;\theta)$ can be obtained by using the standard Itô approach as well. All other cases are out of the reach of the standard theory, as explained in Sec. 5.

**Figs. 6-10** show the evolution of the one-time response moment functions $m_x(t)$, $C_{xy}(t,t)$ and $C_{xx}(t,t) = \sigma_x^2(t)$, obtained both by the numerical solution of the causal system (17), (27) (solid lines ——) and by the MC simulations (dashed lines - - -). In all cases discussed herewith the initial values of the considered moments have been taken to be $m_{x_0} = 2$, $C_{x_0 y_0} = 0$ and $C_{x_0 x_0} = 1$. The three cases shown in each (out of six) subplot of each figure, correspond to three different values of an important parameter: the nonlinearity parameter $\mu_3$, in **Figs. 6** and **7**, the correlation time $\tau_{yy}^{\text{corr}}$, in **Fig. 8**, the central spectral frequency $\omega_0$, in **Fig. 9**, and the parameter $\kappa_3$, controlling the non-Gaussian character of the excitation, in **Fig. 10**. The three cases are distinguished by three different markers (○, △, □). The headline of each figure, in conjunction with the in-figure legend, provide a complete description of data values, fully identifying each case. The filled markers (●, ▲, ■), appearing in the right vertical axis of each subplot, indicate the corresponding long-time statistical equilibrium limits $C_{xy}^{(\infty)}(t,t) = C_{xy}^{(\infty)}(0)$ and $C_{xx}(t,t) = C_{xx}^{(\infty)}(0)$, as obtained by solving the moment system (17), (27) directly in the long time, by means of an analytical solution method which will be presented in a forthcoming paper [68]. These values are always in excellent agreement with the corresponding numerical long-time results, obtained by solving the moment system.

In **Fig. 6** the effect of the (monostable) nonlinearity is investigated, under Gaussian excitation ($\kappa_3 = 0$). In the linear case ($\mu_3 = 0$) the results from the moment system are in full agreement with MC results, as well they should. The agreement is very good in the nonlinear cases as well. Stronger nonlinearity leads to smaller values of the response moments. This fact is expected, since the nonlinear term has negative coefficient ($\mu_3 < 0$), amplifying the restoring effect of the linear term (monostability). Similar results are shown in **Fig. 7**, for some cases of non-Gaussian excitations ($\kappa_3 = 0.4$). The trends are much the same with those appearing in the previous figure, but now the agreement between the moment system solutions and MC results is poorer, because of the non-Gaussian character of the excitation. In the linear case ($\mu_3 = 0$), the non-Gaussian effect is stronger and is suppressed as the nonlinearity increases, improving the agreement with the MC results. Comparing the results presented in **Figs. 6** and **7**, we may conclude that, for a monostable cubic half-oscillator, the effect of the input non-Gaussianity is stronger than the effect of the nonlinearity, with respect to the Gaussian closure adequacy.

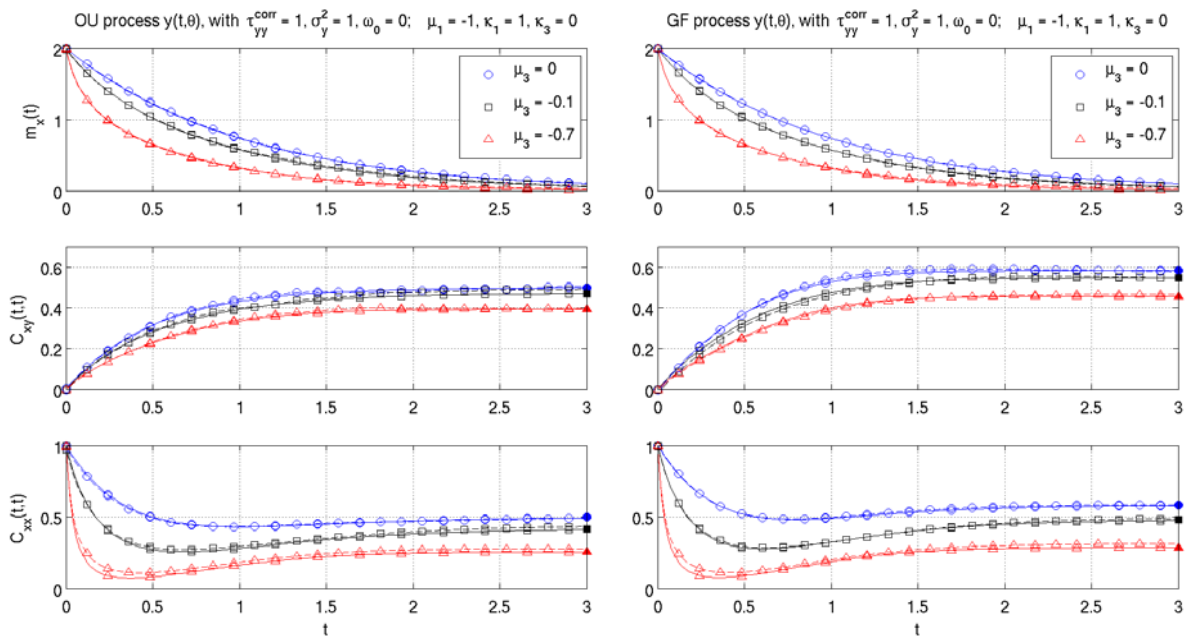

**Fig. 6.** Evolution of one-time response moments $m_x(t)$, $C_{xy}(t,t)$ and $C_{xx}(t,t)$, for a cubic half-oscillator under Gaussian excitation ($\kappa_3 = 0$), for three values of the nonlinearity parameter $\mu_3$. Solid lines (———) indicate results obtained by solving the moment system, and dashed lines (– – –) indicate results obtained by MC simulations. The excitation $z(t;\theta) = \kappa_1 \cdot y(t;\theta)$ is an OU process (left panel) or a Gf process (right panel). The filled markers (●, ▲, ■), appearing in the right vertical axis of each plot, indicate the corresponding long-time statistical equilibrium limits $C_{xy}^{(\infty)}(t,t) = C_{xy}^{(\infty)}(0)$ and $C_{xx}^{(\infty)}(t,t) = C_{xx}^{(\infty)}(0)$, as obtained by solving analytically the moment system in the long time.

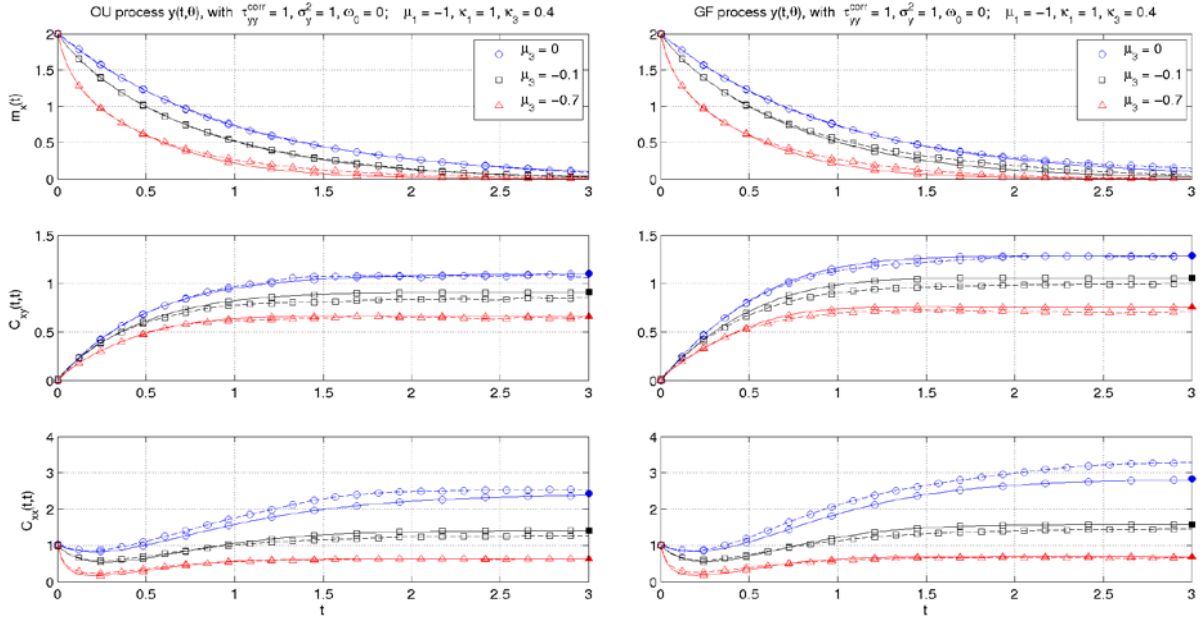

**Fig. 7.** The same as Fig. 6, except that the excitation $z(t;\theta) = \kappa_1 \cdot y(t;\theta) + \kappa_3 \cdot y^3(t;\theta)$ is now non-Gaussian ($\kappa_3 = 0.4$).

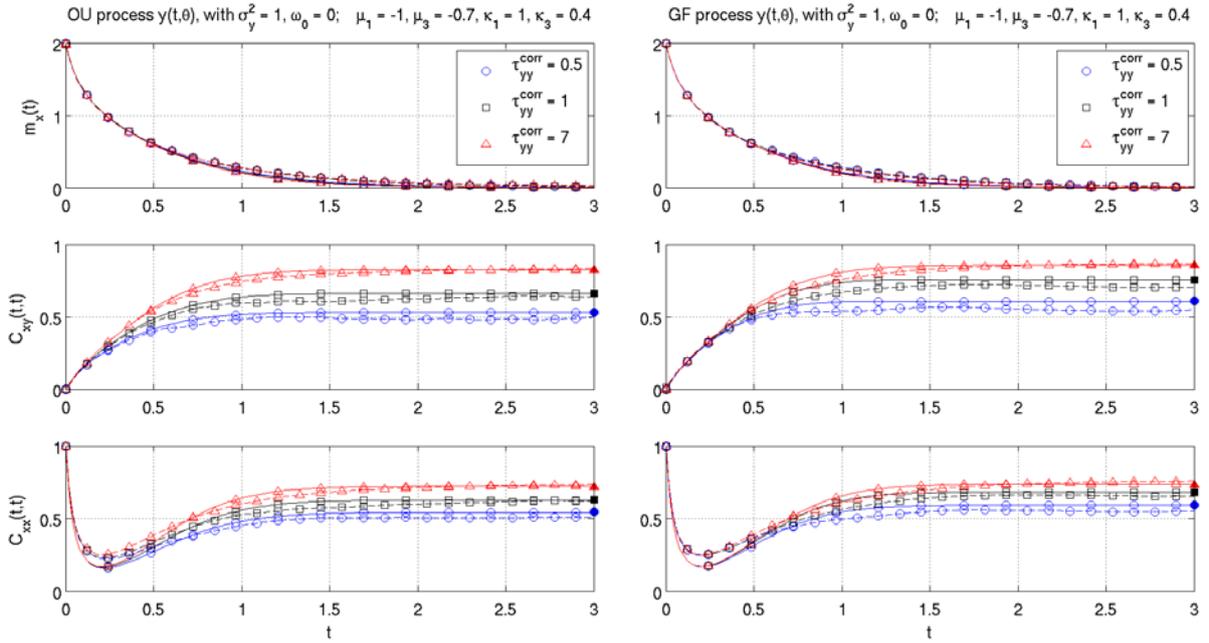

**Fig. 8.** Evolution of one-time response moments $m_x(t)$, $C_{xy}(t,t)$ and $C_{xx}(t,t)$, for a cubic half-oscillator under non-Gaussian excitation ($\kappa_3 = 0.4$), for three values of the correlation time $\tau_{yy}^{corr}$. The excitation $z(t;\theta) = \kappa_1 \cdot y(t;\theta) + \kappa_3 \cdot y^3(t;\theta)$ is derived by an OU process $y(t;\theta)$ (left panel) or by a Gf process (right panel). All conventions regarding the style of the lines and the meaning of the markers in the plots are the same as in **Fig. 6**.

In **Fig. 8** the effect of the correlation time $\tau_{yy}^{corr}$ on the one-time moments is investigated. Mean value $m_x(t)$ is quite insensitive, in contrast with 2$^{nd}$-order moments which exhibit appreciable variation as $\tau_{yy}^{corr}$ changes. The larger the correlation time the larger the values of the moments $C_{xy}(t,t)$ $C_{xx}(t,t)$; see also **Fig. 12**. The differences between the results obtained by the moment system and by MC simulations are negligible for large values of the correlation time, always remaining quite acceptable. The larger differences, observed for the smaller correlation time, are about 7%. These differences reflect the degree of inadequacy of the Gaussian closure, in the presence of nonlinearity ($\mu_3 = -0.7$) and non-Gaussianity of the excitation ($\kappa_3 = 0.4$).

In **Fig. 9** the evolution of the moments $m_x(t)$, $C_{xy}(t,t)$ and $C_{xx}(t,t)$ is shown for three different values of the central spectral frequency $\omega_0$ of the input random function $y(t;\theta)$. Again the mean value $m_x(t)$ is not sensitive. Both 2$^{nd}$-order moments $C_{xy}(t,t)$ and $C_{xx}(t,t)$ decrease as the central frequency $\omega_0$ increases. The agreement between the moment system solutions and the MC simulations is very satisfactory in both Gf and OU cases. Similar results are shown in **Fig. 10**, for three different values of parameter $\kappa_3$. Both $C_{xy}(t,t)$ and $C_{xx}(t,t)$ decrease monotonically with $\kappa_3$, becoming nearly zero for the strongly non-Gaussian case $\kappa_3 = -0.4$, indicated in the figures by ─○─ . The most important finding here is the total failure of the moment system to represent correctly the response moments in the case $\kappa_3 = -0.4$. Clearly, this should be attributed to the Gaussian closure; however, the discrepancy in this case is much stronger than in all other cases. This fact calls for a more detailed analysis. In order to understand the reason for such a strong discrepancy, we have performed a direct check of the quality of the Gaussian closure, by estimating the joint, RE pdf $f_{x(t)\,y(t)}(x,y)$ (in the long time), for the three cases shown to the left panel of **Fig. 10**, where $y(t;\theta)$ is an OU process. The corresponding results are presented in **Fig. 11**. Observing the three subplots of this figure, we easily identify the essentially different structure of the joint RE pdf $f_{x(t)\,y(t)}(x,y)$ in the case $\kappa_3 = -0.4$, in comparison with the other two cases. The two heavy tails that are developed in the case $\kappa_3 = -0.4$, in a direction almost orthogonal to $x=y$, indicate that the Gaussian closure should be inappropriate for this case. To quantify this assertion we have calculated the long-time limiting values of the 4$^{th}$-order moments $C_{xy}^{13}(t,t)$, $C_{xy}^{31}(t,t)$ and the 2$^{nd}$-order moments $C_{yy}(t,t)$, $C_{xy}(t,t)$, $C_{xx}(t,t)$ directly from the MC pdf, and formed the ratios $r_{xy}^{13} = C_{xy}^{13}(t,t)/\left(3 \cdot C_{yy}(t,t) \cdot C_{xy}(t,t)\right)$ and $r_{xy}^{31} = C_{xy}^{31}(t,t)/\left(3 \cdot C_{xx}(t,t) \cdot C_{xy}(t,t)\right)$, which are assumed equal to 1 under Gaussian closure. The calculated values of these ratios are shown in **Table 1**. For the cases $\kappa_3 = +0.4$ and $\kappa_3 = 0$ the values of the ratios are close to unity, justifying the Gaussian closure assumption. For the case $\kappa_3 = -0.4$ the values of $r_{xy}^{13}$ and $r_{xy}^{31}$ exhibit a huge departure from unity, explaining the failure of the Gaussian closure.

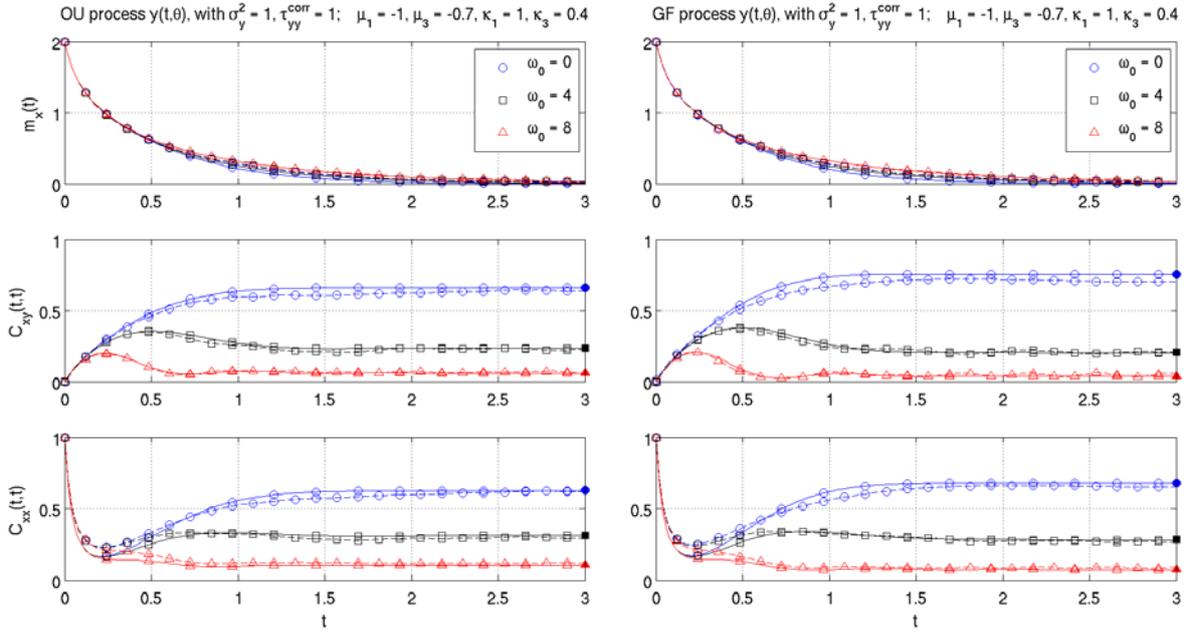

**Fig. 9.** The same as **Fig. 8**, except that now the varying parameter is the central spectral frequency $\omega_0$.

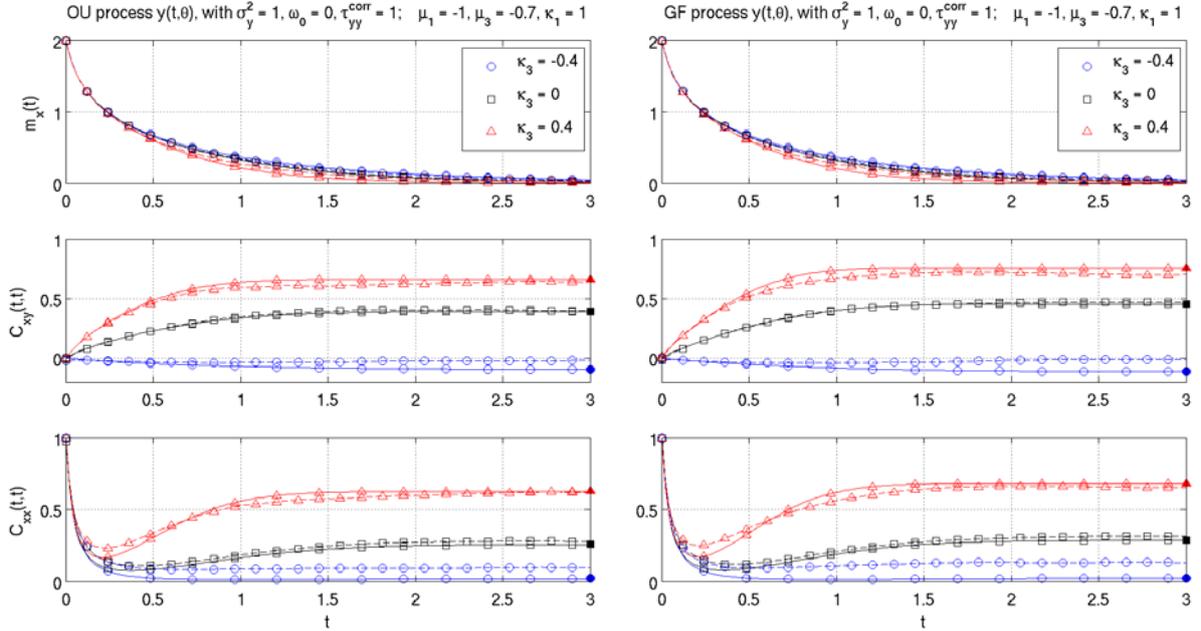

**Fig. 10.** The same as **Fig. 8**, except that now the varying parameter is $\kappa_3$.

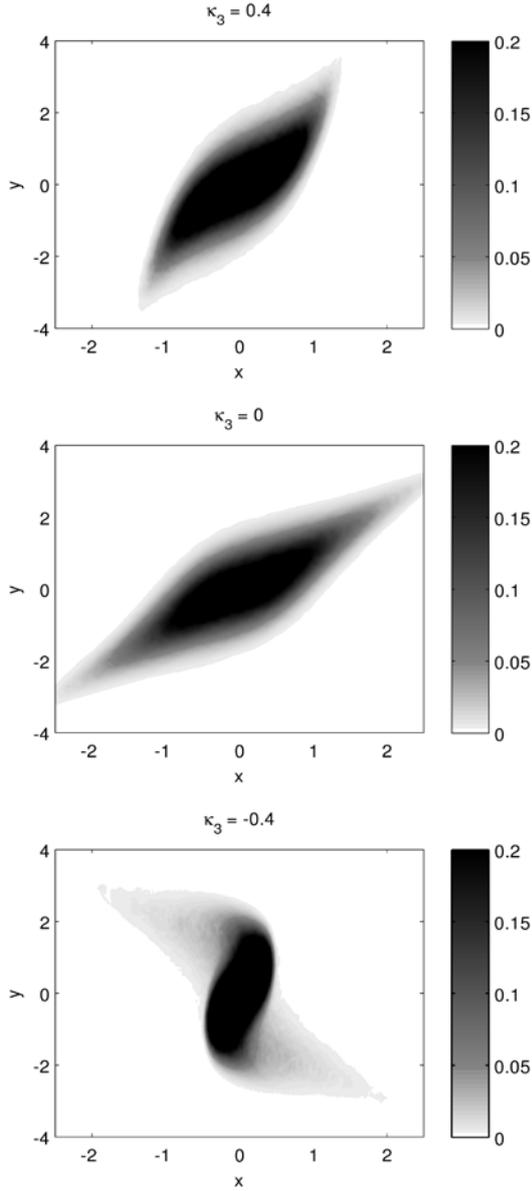

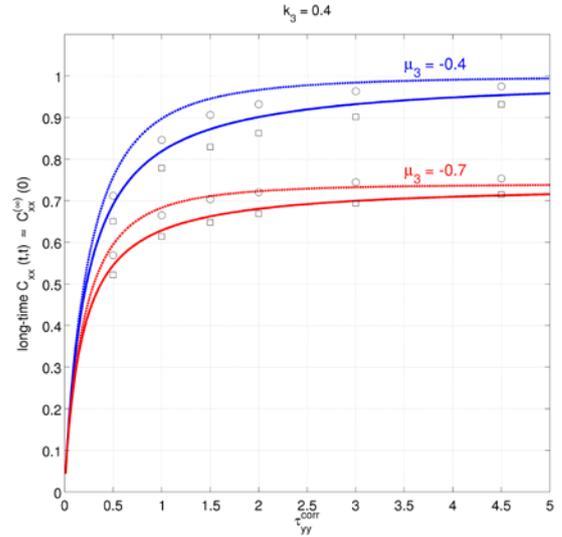

**Table 1**

| $\kappa_3$ | $r_{xy}^{13}$ | $r_{xy}^{31}$ |
|---|---|---|
| $+0.4$ | 1.08 | 1.10 |
| 0 | 0.96 | 0.78 |
| $-0.4$ | 7.45 | 15.85 |

**Fig. 11.** The long-time, joint RE pdf $f_{x(t)\,y(t)}(x, y)$, for the three cases shown in the left panel of **Fig. 9**, calculated by MC simulations.

**Fig. 12.** Long time variance $C_{xx}(t,t) \approx C_{xx}^{(\infty)}(0) = \sigma_x^2(\infty)$ for a cubic half-oscillator under Gf (dashed lines) and OU (solid lines) non-Gaussian input ($\kappa_3 = 0.4$), with $\sigma_y^2 = 1$ and $\omega_0 = 0$, versus $\tau_{yy}^{\text{corr}}$. Results for two values of the nonlinearity parameter ($\mu_3 = -0.4$ and $\mu_3 = -0.7$) are shown. Corresponding MC simulation results are also depicted for Gf (circle markers) and OU (square markers) random input.

In **Fig. 12** the effect of the variation of the shape of the covariance of the input random function $y(t;\theta)$ (Gf vs OU with the same defining parameters $\sigma_y^2$, $\tau_{yy}^{\text{corr}}$, $\omega_0$) to the long-time response variance, $\sigma_x^2(\infty) = C_{xx}^{(\infty)}(0) \approx C_{xx}(t,t)$ for values of $t$ in the long-time regime, is illustrated. The plotted variances are obtained by the long-time numerical results of the moment system under Gf and OU input, with $\sigma_y^2 = 1$, $\omega_0 = 0$ and $0 < \tau_{yy}^{\text{corr}} < 5$. Two cases are plotted, corresponding to a mildly nonlinear ($\mu_3 = -0.4$) and a strongly nonlinear ($\mu_3 = -0.7$) monostable, half oscillator, under unimodal, non-Gaussian excitation ($\kappa_3 = -0.4$). The moment system solutions are compared with MC simulations, also shown in **Fig. 12**. The two methods yield similar results, with relative differences between 1% and 7%. In agreement with findings already discussed previously, the long-time response variance $C_{xx}^{(\infty)}(0)$ decreases with the nonlinearity (see also **Fig. 7**), and increases with the correlation time $\tau_{yy}^{\text{corr}}$ of the input random function (see also **Fig. 8**). The most important finding here is that the long-time response variance $C_{xx}^{(\infty)}(0)$ is significantly affected by the shape of the input function (Gf vs OU), especially in the vicinity of $\tau_{yy}^{\text{corr}} = 1$, despite the fact that all integral parameters ($\sigma_y^2$, $\tau_{yy}^{\text{corr}}$, $\omega_0$) are the same. The response variance obtained under Gf random input is always higher than the response variance obtained under OU input. Around $\tau_{yy}^{\text{corr}} = 1$ the difference is as high as 10%. This finding highlights the importance of methods treating random ODEs under general random excitation, beyond the limitations of Itô/FPK approach. Recall that only the case of a centered OU input can be modeled and solved within the standard Itô/FPK approach (Sec. 5).

In **Fig. 13**, the two-time moments $C_{xy}(t,s)$ and $C_{xx}(t,s)$ are plotted as three-dimensional surfaces over the whole two-time rectangle $R_{ts}(T) = \{(t,s) : t_0 \leq t \leq T, t_0 \leq s \leq T\}$, for the case $\mu_1 = -1.0$, $\mu_3 = -0.7$, $\kappa_1 = 1.0$, $\kappa_3 = 0.4$, with a Gf input function $y(t;\theta)$ defined by $\sigma_y^2 = 1.0$, $\omega_0 = 0.0$, $\tau_{yy}^{\text{corr}} = 1.0$. Results obtained both by numerical solution to the moment system (left panels) and by MC simulations (right panels) are shown. Qualitatively, the corresponding surfaces are very similar. For quantitative comparisons, intersections of the surfaces with planes parallel to the $t$-axis are plotted (see **Fig. 14**) and discussed in the next paragraph. The one-time moments $C_{xy}(t,t)$ and $C_{xx}(t,t)$ are obtained from this figure as the intersection of the surfaces with vertical planes passing through the line $t = s$ of the base rectangle, and they are shown by dashed lines. The stationary, long-time moments $C_{xy}^{(\infty)}(t-s)$ and $C_{xx}^{(\infty)}(t-s)$ are shown by solid curves in the left panel of **Fig. 13**, as have been obtained by an asymptotic analysis of the moment system, which will be presented elsewhere [68]. It is clear that they nicely fit to the numerical solution obtained by solving the moment system.

Finally, in **Fig. 14**, the two-time moments $C_{xy}(t,s)$ (left panel) and $C_{xx}(t,s)$ (right panel) are plotted as functions of the response time $t$, for various values of the excitation time $s$. Results obtained by the moment system (solid lined) and by MC simulations (dashed lines) are generally in good agreement, with values near the peaks slightly overestimated by the moment system, the mismatch varying between 3% and 10%. For $s=3$ $(=3\times\tau_{yy}^{corr})$, the system is no longer affected by the initial conditions, having reached the long-time, statistical equilibrium state. This is why the response covariance $C_{xx}(t,s)$ (right panel) is symmetric with respect to the line $t=s=3$, as expected. In all cases the cross-covariance $C_{xy}(t,s)$ (left panel) reaches a peak for $t>s$, and has steeper decay in the side $t<s$. Note that, for $s>0$, $C_{xy}(t,s)$ has non-zero values for a significant time period when the excitation is in advance of the response ($t<s$). This is due to the nature of the random input considered. Its covariance being smooth, leads to a response smoothly correlated both with past and future values of the excitation. As a result $C_{xy}(t,s)$ is non-zero for a time period proportional to $\tau_{yy}^{corr}$. In the examined cases, where $\tau_{yy}^{corr}=1.0$, the response stays correlated with the input for about $2\times\tau_{yy}^{corr}$ ahead ($t<s$) and $3\times\tau_{yy}^{corr}$ behind ($t>s$).

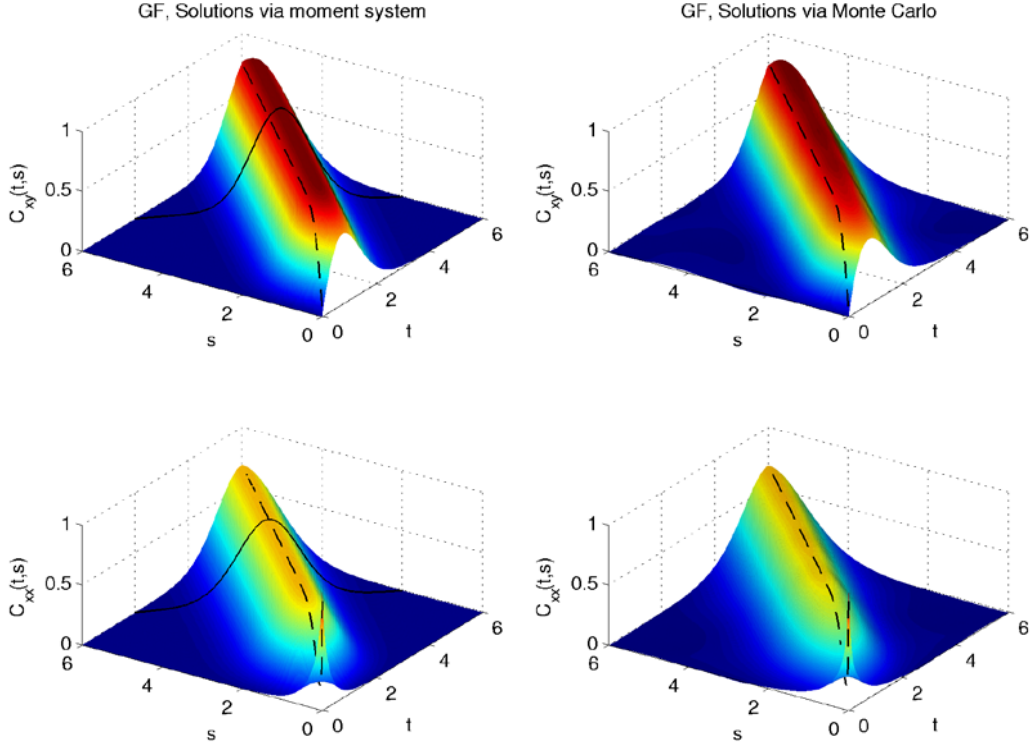

**Fig. 13.** Two-time moments $C_{xy}(t,s)$ and $C_{xx}(t,s)$ calculated by numerical solution to the moment system (left panel) and by MC simulations (right panel). Dashed curves represent the intersection of the surfaces with vertical planes passing through the line $t=s$ and correspond to the one-time moments $C_{xy}(t,t)$ and $C_{xx}(t,t)$. Solid curves correspond to the stationary long-time moments $C_{xy}^{(\infty)}(t-s)$ (upper left panel) and $C_{xx}^{(\infty)}(t-s)$ (lower left panel), obtained analytically.

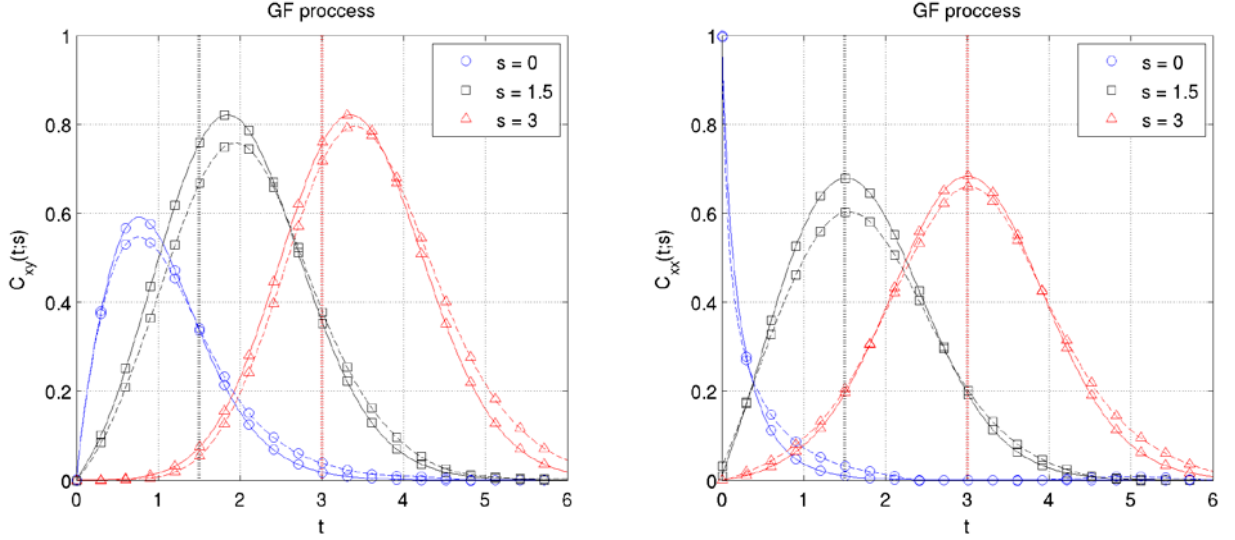

**Fig. 14.** Evolution of the two-time moments $C_{xy}(t,s)$ (left panel) and $C_{xx}(t,s)$ (right panel) with respect to the response time $t$, for three values of the excitation time $s=0$, $s=1.5$, $s=3$. The curves represent intersections of the surfaces of Fig. 13 with vertical planes parallel to the $t-$axis. Solid curves indicate results obtained by solving the moment system and dashed curves indicate results obtained by MC simulations. The dashed vertical lines indicate the state $t=s$.

## 7. Discussion and conclusions

In this paper we have introduced the two-time, response-excitation moment equations for a cubic half-oscillator, excited by colored (Gaussian or non-Gaussian) noise, and demonstrated that these equations can be recast as a closed, causal (non-local in time), solvable system. To achieve this status a two-fold closure is necessary: moment closure and time closure. Moment closure is approximate and leans on well-known ideas and techniques. Time closure is exact (given the preceding moment closure), and is made possible by appropriate exploitation of the inherent structure of the moment system. Notwithstanding the specific application of the present work, the methodology described herewith is quite general; it can be applied to full oscillators and other higher-order dynamical systems.

The proposed methodology, of formulating moment equations for non-linear systems under smoothly-correlated excitation, generalizes the usual Itô/filtering approach, in the sense that it applies to excitations having arbitrary correlation structure. When the input random function can be represented as a solution of an Itô SDE the causal (non-local) moment system obtained by the present method simplifies to a local one, becoming identical with the one-time moment system obtained by the Itô/filtering approach.

The causal nonlinear moment system, although unusual in the context of moment equations for dynamical systems, belongs to a family of functional differential equations which has been studied for almost a century, since its introduction by Volterra and Tonelli. The solvability (well-posedness) of the causal moment system is guaranteed by the existing mathematical theory, and its numerical solution, developed in this paper, is both efficient and robust. The two-cycle (multi-scale) approach to numerical solution exhibits an astonishing self-correcting feature, leading to very accurate numerical solutions. Numerical results indicate that the output two-time moments are significantly affected both by the correlation time and the details of the shape of the input random function, which in general cannot be taken into account in the Itô/FPK/filtering approaches. The obtained results compare satisfactorily with extensive Monte Carlo simulations results.

The two-time moment equations presented herewith constitute a building block of the more general response-excitation theory, introduced in 2006 by Athanassoulis and Sapsis [50], [52], and aiming at the calculation of the joint, RE pdf. In fact, this type of moment equations, when applied locally in the RE phase space, provides closure conditions for the differential constraint satisfied by the joint RE pdf, permitting its unique solution. See [55], [56] for first results in this direction, which will be fully elaborated and presented in forthcoming papers.

## Acknowledgments

Ivi C. Tsantili and Zacharias G. Kapelonis are PhD candidates at the School of Naval Architecture and Marine Engineering, National Technical University of Athens (NTUA), supported by NTUA scholarships.